\documentclass{article}
\usepackage{inputenc}

\title{\Large{A Probabilistic Framework for Optimizing Real-time Decisions in Humanitarian Aid Delivery Systems}\footnote{Incomplete work in progress.}}

\author{Roozbeh Yousefzadeh \vspace{.2cm}\\
\vspace{.2cm} Yale University, New Haven, CT\\
\texttt{roozbeh.yousefzadeh@yale.edu}}
\date{}

\usepackage{indentfirst}
\usepackage{graphicx}
\usepackage{titlesec}
\usepackage{bbm}
\usepackage[margin=0.8in]{geometry}
\usepackage{comment}
\usepackage{setspace}
\usepackage{amsmath}
\usepackage{amssymb}
\usepackage{relsize}
\usepackage[usenames, dvipsnames]{xcolor}
\usepackage[skip=0pt]{caption}
\usepackage{float}
\usepackage{tikz}
\usetikzlibrary{positioning}
\usepackage[ruled,vlined,linesnumbered]{algorithm2e}

\usepackage{natbib}
\setcitestyle{authoryear,open={(},close={)}}

\usepackage{mathrsfs} 
\usepackage{hyperref}
\hypersetup{
    backref=true,
    colorlinks=true,
    citecolor=blue,
    citebordercolor=white}
    
\DeclareMathOperator{\erf}{erf}
\usepackage{caption}
\usepackage{subcaption}
\usepackage{cleveref}
\begin{document}
\maketitle
\onehalfspacing
\large

\section*{Abstract}
This paper presents a computationally efficient model for optimizing real-time decisions in humanitarian aid delivery systems. Our formulation models a hierarchical system and is a mixed integer, probabilistic, non-linear and non-concave optimization problem. The proposed model considers the costs and probabilistic nature of transfer times and maximizes the reliability of the system using the available budget. We account for late deliveries using a nonlinear penalty function. We also propose an algorithm that uses a directed acyclic graph to deal with the discrete variables in tandem with a homotopy method for optimizing the continuous variables. We then offer a pruning method to eliminate the cost inefficiencies in the system. The effectiveness of formulation is examined in some numerical examples in development.

\section{Introduction} \label{introduction}

Humanitarian aid delivery has been the subject of research, especially due to the rise of humanitarian crises around the world. Generally many different agencies and organizations are involved in delivering different kinds of assistance and aid to the people in need in an on-going basis. One of the main aspects of this problem is that time is of the essence. In other words the main goal is to ensure that sufficient aid is delivered at the time it is needed. Obviously late delivery of aid can have devastating humanitarian consequences. 

In order to achieve this goal, it is necessary to have a prepared network of aid delivery in place. We refer to this aspect as \textit{pre-disaster planning} and it involves the design of the system (e.g., location and characteristics of facilities and vehicles, shape of the network including the hubs and spokes) and the pre-positioning of aid commodities within the network.

The other prominent aspect of the humanitarian aid delivery concerns operations that should happen in the aftermath of disasters. We refer to this aspect as the \textit{post-disaster operations}, the process of transporting the aid commodities through the system and delivering it to the disaster areas where there is demand for aid. These operations have to be performed as fast as possible and the operational decisions have to be made in real-time. These real-time decisions are mainly about the amount of aid commodities dispatched with the available vehicles, and how the vehicles are routed through the network.

In most instances, the military is responsible for these operations with coordination with other organizations such as the Red Cross. The United States Department of Defense has specific programs and infrastructure to conduct such delivery operations, which is referred to as HA/DR: the Humanitarian Assistance and Disaster Relief Program (\cite{o2012foreign,tarnoff2009foreign}). 

While the system is in operation, there are many uncertainties mainly with respect to the completion time of transportation activities. The actual time it takes for a plane, for example, to travel between two points is often different from the expected travel time due to many possible reasons such as weather. Nevertheless, these uncertainties can be modeled relatively accurate by using probability distributions.

Given all the operational uncertainties and the complexities in the network (\cite{van2006humanitarian}), one has to rely on a mathematical model in order to make the real-time decisions optimal. But such a mathematical model should be computationally efficient so that it can be solved in real-time.

An aid delivery system can rely on many modes of transportation. Among different modes, aerial delivery is one of the most effective ways of delivering aid after severe disasters in hard to reach regions. In this method of aid delivery, planes fly at low altitude so that the pilot can visually inspect the situation on the ground and drop the aid with some precision.

Based on the existing infrastructure and available fleet of vehicles in the affected regions, there might be other modes of transportation available to deliver the aid commodities as well (e.g. trains, trucks, drones, etc). Therefore in our study we consider a general inter-modal transportation system for the delivery of aid commodities. We focus on the problem of real-time decisions regarding the flow of aid in a hierarchical network.

In section \ref{literature} we review the literature. In section \ref{problem} we describe the problem and different aspects of real-time decision making in an aid delivery system. In section \ref{model} we describe the model and present the formulation. Section \ref{results} contains the results obtained from the model and section \ref{conclusion} is the conclusion.


\section{Literature review} \label{literature}

Here we review the literature on humanitarian aid delivery systems with a focus on methodologies related to optimizing the real-time decisions.

\cite{haghani1996formulation} formulated a multi-commodity multi-modal network flow model for disaster relief operations. They considered the flow in the network during time windows and proposed two heuristic algorithms for solving their proposed formulation. The objective in their study is minimizing the total cost, which includes the vehicular and commodity flow costs, the supply carry-over costs, and the transfer costs. Their formulation is a mixed integer linear program (MILP) which does not consider the probabilistic nature of transfer times.

\cite{barbarosouglu2002interactive} proposed a mathematical model for helicopter mission planning during a disaster relief operation. Their formulation consists of two mixed integer programs, one dealing with fleet assignments and tour numbers for each helicopter and the other dealing with the vehicle routing and loading decisions. They optimize the sub-problems iteratively. However, their network is not hierarchical and their transfer times are not probabilistic.

\cite{barbarosoglu2004two} extended the work by \cite{haghani1996formulation} to consider probabilistic demand for aid at the destinations. The objective in their study is to minimize the total cost, and the transfer times in the network are not considered to be probabilistic.

\cite{sheu2007emergency} studied the problem of emergency disaster relief operations and focused on the quickness of aid delivery. The proposed model has a probabilistic approach in quantifying the demand but does not consider the probabilistic transfer times in the network. The network itself has one layer of hubs between the origins and destinations, rather than the general hierarchical network in our model.

\cite{vitoriano} proposed a multi-criteria optimization model which considers not only the total cost in the system but also takes into account the time of response, equity of distribution among destinations, and security of operation routes. They use goal programming techniques for their formulation and only consider a single aid commodity. Vehicles in the model go through several links to reach their destinations but vehicles do not interact and the network is not hierarchical. They also mention probabilities of successful transfer for each link in the network but the probabilities are perceptive and not calculated mathematically.



\cite{huang2012models} considered objectives beyond minimizing the cost. They defined and formulated performance metrics representing efficacy (i.e., the extent to which the goals of quick and sufficient distribution are met), equity (i.e., the extent to which all recipients receive comparable service) and efficiency (i.e., minimizing the cost). Their formulation is a deterministic integer program.

\cite{afshar2012modeling} proposed a model that controls the flow of several relief commodities in an aid delivery system. The objective in their study is to minimize the total amount of unsatisfied demand in the system. Their formulation is a deterministic linear integer program.

\cite{zary} investigated the literature on humanitarian logistics published between 2001 and 2014 and performed a bibliometric analysis regarding the citations and co-citations of the publications. The goal of this study was to identify the knowledge network among this field of research and to summarize relevant theories, concepts and research methods.


\cite{das2014relief} studied the problem of inventory management (pre-positioning of aid in warehouses) for unforeseen humanitarian disaster relief operations. They proposed a stochastic model that optimizes the amount of inventory that should be kept in a warehouse to satisfy a probabilistic demand with probabilistic lead-time. They use uniform probability distributions and do not consider any capacity constraints. Although their research tackles a planning aspect of humanitarian aid delivery, its model could use our real-time model to better estimate how long it would take for the aid to reach its destination in a hierarchical network.

\cite{ozdamar2015models} reviewed the literature on the response and recovery planning phases of the disaster relief operations. They classified the mathematical models in terms of vehicle/network representation structures and their functionality. It can be seen that the inter-modal systems of delivery in humanitarian crisis have rarely been studied.

\cite{huang2015modeling} considered multiple humanitarian objectives in emergency response to disasters. Their objectives are life saving utility, delay cost and fairness, and a time-space network is used to optimize the system. The approach to put the humanitarian objectives front and center is novel. However, they neglect probabilities and their network is not hierarchical.

\cite{lu2016real} proposed a rolling horizon-based framework for real-time relief distribution in the aftermath of disasters. The objective in this study is to minimize the total time to deliver the aid to satisfy the demand. The formulation is a linear integer program, the network is not hierarchical and they neglect the probabilities in the system.

\cite{rezaei2016interactive} studied the problem of inventory management with respect to perishable aid commodities (e.g. medicine). They considered a lifetime for such commodities and developed a mixed integer linear program to optimize the amount of perishable aid to be kept in the warehouses and to find the best policy for renewing the stocked commodities at the pre-disaster phase. Their objective is to minimize both the average response time and the costs in the system. They consider probabilities to estimate the chance of a disaster occurrence but their study is limited to inventory management inside a warehouse.

\cite{nagurney2016generalized} developed a generalized Nash equilibrium network model for post-disaster relief operations. This study uses game theory to model the competition between multiple non governmental organizations (NGOs) over financial funds in response to disasters. Different disaster relief strategies are investigated from the view point of the coordinating authorities (such as the United Nations) and the NGOs.

\cite{balcik2016literature} reviewed the literature focused on inventory management, facility location and planning aspects of humanitarian aid delivery which is insightful but not directly applicable to the problem of real-time decisions.

Some studies such as \cite{chiu2007real}, \cite{huang2013continuous} and \cite{balcik2017site} focus on real-time routing problem for ground modes of transportation.
\cite{huang2013continuous} for example, proposes an integer program with the objective to minimize the sum of arrival times. However, since our study is about inter-modal aid delivery systems, the routing can be seen as a separable problem. In other word, the best route found by the routing models can be used as the input for our model.

\cite{army,bastian} specifically studied the problem of aerial aid delivery by the military. They formulated a mixed integer program and used goal programming to optimize the system with respect to three goals: target response time, target budget, and target demand. The objective in their study is to minimize the sum of deviations from the goals. The network in their study is hierarchical, and they consider probabilistic demand since their focus is mainly on pre-disaster planning aspects of the problem. They neglect the probabilistic transfer times and the formulation is only applicable for a single aid commodity system.

\cite{tardiff2017development} designed an airdrop technology to successfully deliver humanitarian aid over populated areas. That study considered injury thresholds of free-falling aid items and various types of aid that needs to be included in a drop package. They also developed practical guidelines for aerial aid delivery and a cost assessment tool.

The problem of optimizing real-time decisions in freight systems has been the subject of recent study in contexts other than aid delivery, such as for ship liners by \cite{li2016real}, for high speed trains by \cite{zhan2015real}, for railways by \cite{dollevoet2017application}, for high frequency transit by \cite{sanchez2016real}, for express urban pick-up and delivery by \cite{ferrucci2014real} and for buses by \cite{berrebi2015real}. Although all these studies are insightful about optimizing real-time decisions in transportation systems, their objectives are different from objectives pursued for the humanitarian aid delivery.

Recently, \cite{yousefzadeh2019probabilistic} developed a homotopy method and a probabilistic framework to optimize real-time dispatch decisions in hierarchical networks. The problem they consider, however, is different than our problem here, because they only consider optimizing the dispatch times, while we consider dispatch times in tandem with other discrete variables.\footnote{The way in which a typical freight system is operated is also different than the way an aid delivery system is operated. In aid delivery the goal is to mobilize the system as efficiently as possible, given many resource constraints, and try to satisfy the demand at downstream as much as possible. But, in a typical freight system, we are dealing more or less with a network flow problem in which the flow is initialized by customers at the upstream, and the goal is to minimize the total cost of delivering.}

In summary, researchers have studied many different aspects of aid delivery systems. However, our work is novel in investigating real-time decisions in a hierarchical network considering the probabilistic nature of transfer times and the imminent need for reliable delivery of aid.

\section{Problem Statement} \label{problem}

In this paper we consider the United States humanitarian assistance and disaster relief (HA/DR) operations as previously studied by \cite{army,bastian} and focus on optimizing the real-time decisions while considering probabilistic transfer times in an inter-modal transportation system.

\begin{figure}[h]
\centering
\includegraphics[scale=0.7]{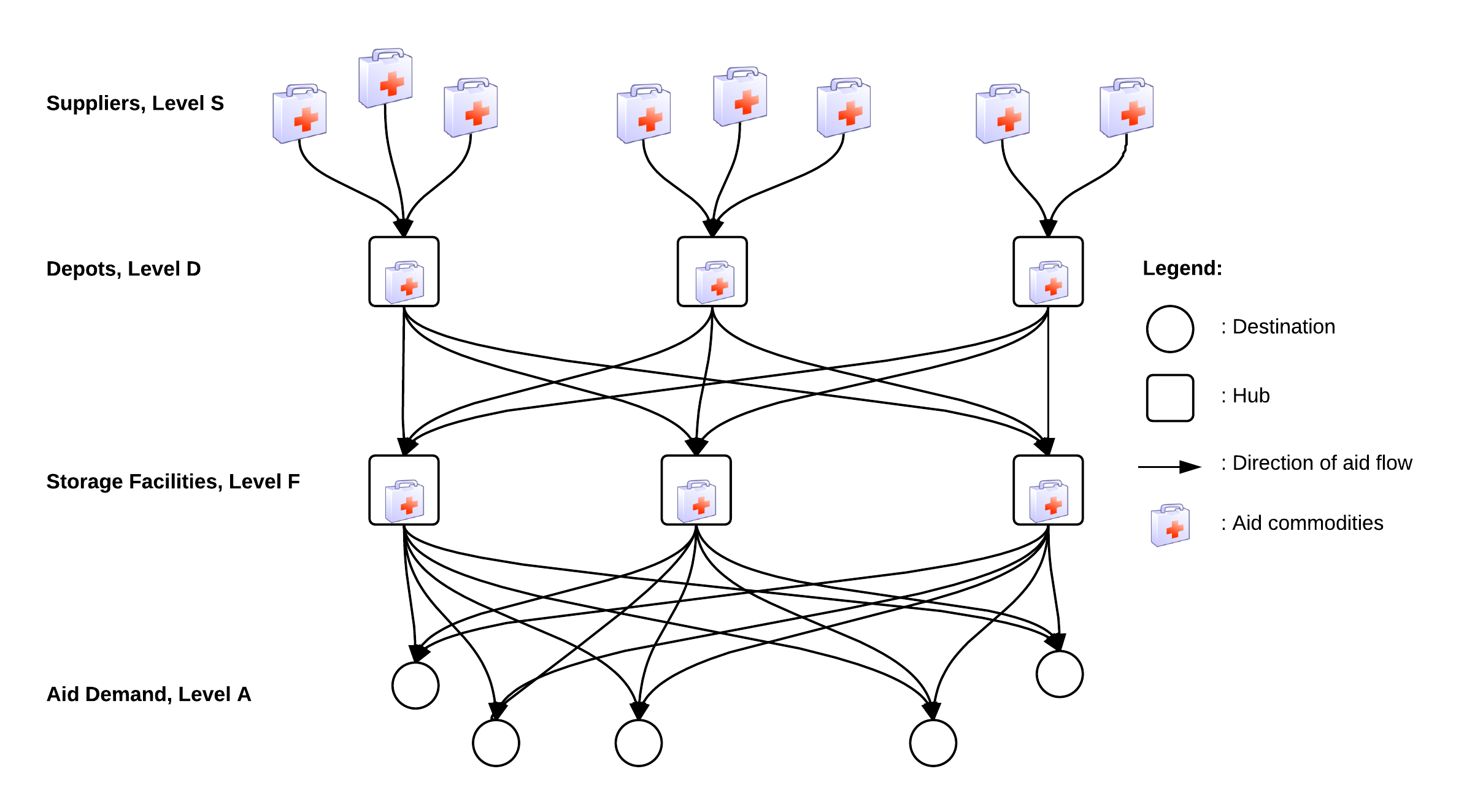}
\caption{Sketch of a generic hub-and-spoke aid delivery system} 
\label{fig:network}
\end{figure}

As an example, consider the network in Figure~\ref{fig:network} and aim to optimize the decisions that have to be made when responding to the demand at the destinations. The destinations on level $A$ initiate the flow of aid commodities in the network. For \textit{pre-disaster planning} purposes, the amount of aid commodities that are pre-positioned through the network can be determined based on probabilistic demand in the future. But in the \textit{post-disaster phase}, the actual delivery of commodities at the destinations should be based on verified imminent need for aid. Therefore in our real-time problem, demand is described as the certain quantity of aid commodities that should be delivered within a given time frame. 

The storage facilities on level $F$ are where the delivery vehicles are located. The delivery vehicles can belong to any transportation mode. The time required to perform the transport activities throughout the network will be considered probabilistic and each activity will have its own probability distribution for completion time.

In a real-time instance of the network, some known amount of aid commodities are stored in the storage facilities (level $F$) and depots (level $D$). The aid commodities stored in the depots can be transported to the storage facilities, and additional aid can also be ordered from various suppliers (level $S$) to be sent out to each of the depots. We can classify the real-time decisions based on their location in the network: 

\begin{enumerate}
\item Decisions on level $S$: how much additional aid should be ordered from which suppliers to be delivered to which depots.
\item Decisions on level $D$: which vehicles to dispatch from depots at what time towards which storage facility (level $F$) with what cargo. 
\item Decisions on level $F$: which vehicles to dispatch from storage facilities at what time to what destination and with what cargo. 
\end{enumerate}

As noted in the literature review, the operation of an aid delivery system has different goals that are often in conflict. For example, minimizing the delay for deliveries is likely to have an adverse effect on minimizing the costs. Therefore it is essential to define an effective objective for a model that optimizes the real-time decisions in an aid delivery system.

A common approach that has been pursued in many previous studies is to minimize the total cost in the system and to consider any deviation from the demand as a penalty added to the total cost function. This approach can be effective in commercial freight deliveries because every possible outcome can be monetized as a cost or benefit. However, in humanitarian disasters it will be extremely hard to associate monetary value to late delivery of aid because of the severe consequences. Also, in commercial freight systems, minimizing the operational cost is often a proxy for maximizing the profit, whereas profit is not a concern in disaster relief. 

Another approach for example can be to formulate the total amount of lateness among all deliveries in the system and try to minimize it. The issue with this type of approach is that if some portion of the demand is unsatisfied, the total delay in the system will blow up the objective and make it meaningless.

In practice, often a fixed budget is allocated for aid delivery. Furthermore, unlike commercial freight systems, it is common for an aid delivery system to not be able to satisfy all the demand. Therefore the decision makers have to operate the system such that the best possible outcome is achieved considering the budget and all other limitations in the system (e.g., limitations imposed by the available fleet of vehicles).

We formulate our model based on this practical assumption that the aid delivery system has an allocated budget and the objective is to respond to the demand in the most effective way. In other words the objective is to maximize the reliability of the system with respect to the demand. The reliability can be defined as the weighted probability of on-time delivery of the needed aid commodities. The unit for the weighted probability of on-time delivery is the weight of aid commodities and the ceiling for it is the sum of all the demand in the system. In an ideal case where the system has enough resources to satisfy all the demands with probability 1, the maximal reliability will be achieved. In order to make the reliability measure more meaningful, we normalize the weighted probability of on-time deliveries by the total weight of the demand and call it $\mathcal{R}$ value. $\mathcal{R}$ value or in other words our reliability measure is a continuous variable between 0 and 100 for any aid delivery system.

Because the amount of demand is directly related (perhaps linearly) to the number of affected people, we can conclude that satisfying more demand (i.e., increasing reliability) translates to helping more people and saving more lives. If we further assume that life of human beings is of equal worth regardless of their location, it becomes obvious that our model implicitly tries to treat all the demand (individual human beings) equally.


\section{The model} \label{model}

Our model considers a network like that in Figure~\ref{fig:network} and optimizes the dispatch decisions to maximize the reliability of the system with respect to the demand for aid. The output of the model will be the optimal real-time decisions and the corresponding reliability of the system. The demand at a node is characterized by the amount of each aid commodity that is needed and a time frame for its delivery. The fleet of vehicles available at a node have a time availability which might not be the current time and each vehicle has a specific capacity to carry commodities.

To demonstrate the trade-offs in the system, consider a node $a_1$ on level $A$ with some urgent demand that includes different aid commodities. Suppose, for example, that a node $f_1$ on level $F$ has a vehicle $v^f_1$ available to be dispatched to $a_1$ with enough capacity to carry all the aid needed, but assume there are insufficient commodities available at $f_1$ to satisfy all the demand at $a_1$. Perhaps the shortage can be met from node $d_1$ on level $D$ via vehicle $v^d_1$. In such circumstances we would face the decision of whether to dispatch $v^f_1$ with the aid that is available at $f_1$, or to dispatch $v^d_1$ with the missing commodities and delay the dispatch of $v^f_1$, hoping that $v^d_1$ would arrive with more aid. To complicate this a little, we can imagine that there is another vehicle $v^f_2$ available at node $f_1$ enabling us to satisfy the whole demand with two deliveries: $v^f_1$ delivering the existing commodities at $f_1$ and later $v^f_2$ delivering the aid brought in by $v^d_1$ to node $f_1$. This would obviously increase the operational costs which might make it impossible to satisfy all demands elsewhere in the network.

For any demand, the reliability will be the product of the weight of the dispatched aid and the probability of it arriving on-time. It does not make any sense to deliver excess aid to a destination (on time or late), and therefore the dispatched commodities are either equal to the demand, which is the ideal case, or less than the demand, which adversely affects the reliability. The other factor in the reliability is the probability of on-time delivery. Since the transfer time is considered to have a probability distribution, the probability of on-time delivery will be a function of the length of the time interval between dispatch and delivery.

In the example above, delaying the dispatch of $v^f_1$ is generally not favorable for the probability of on-time arrival at the destination node $a_1$. If there is plenty of time available to satisfy the demand, then delaying the dispatch will have a minimal adverse effect on the probability. However, if the available time frame for $v^f_1$ is tight, then delaying the dispatch can have a very adverse effect on the probability of on-time delivery, and even a short amount of delay in dispatch can plummet the probability to zero. Therefore, when to dispatch $v_1^f$ in our simple example may not be clear without considering a probabilistic model.

If we expand our considerations to include level $S$ and multiple nodes on each of the levels, the situation becomes even more complicated. Then the amount of aid commodities loaded on a single vehicle can possibly be a combination of commodities brought in from various nodes on the upper levels. Any vehicle can also be dispatched to one of various destinations. It is easy to see how the decisions and the trade-offs become more complicated and intertwined as the network grows.

In our model we have formulated this problem for a generalized network with four levels as depicted in Figure~\ref{fig:network}. The main simplifying assumptions in our model are:
\begin{enumerate}
\item Each delivery destination is modeled as a single point. In reality it is plausible that the delivery vehicle will have multiple stops(deliveries) within a small region but our assumption is that small region can be modeled as single node in our network. This assumption is realistic, considering that the capacity of delivery vehicles is usually much smaller than the demand at a single location. Each node on level $A$ of the network should represent a region that can be serviced by one or more delivery vehicles.

\item A vehicle, dispatched towards level $A$ can only serve a single destination node and then will return back to a node on level $F$. The transfer cost accounted for the delivery vehicles is the cost of the round trip.

\item There are no capacity constraints in the storage facilities and depots. Storage capacity is not influential in our model since the model is about real-time decisions and not about the pre-disaster planning aspects of the aid delivery system. In this real-time model, the flow of aid through the network is intended to meet the imminent known demand at the destinations and not intended to fill the storage facilities for unknown demand in the distant future.

\item Both the time and amount of the demand are certain in our network. 

\item The flow of aid in the network is unidirectional and there are no loops.

\item Any probability density functions are applicable to our model as long as they are continuous and differentiable.

\item The probabilistic characteristics of transfer time (e.g. mean and standard deviation) between hubs include the time for transit and off-loading of cargo.

\end{enumerate}

Our formulation is written in tensor notation summarized in \cite{qi2017tensor}. Scalars and vectors (tensors of order zero and one) are denoted by lowercase letters, tensors of order two are denoted by capital letters and tensors of higher order are denoted by Fraktur script letters. Tensors of order greater than zero are denoted by boldface. The relevant level on the network is shown as a superscript in parenthesis. The order of a tensors representing the parameters and variables can vary on different levels of the network. 

Another issue that we have considered in the model is with regards to the probability of late deliveries. We believe the model should distinguish between small and large amounts of lateness. For example delivering the commodities 10 minutes late should not be treated the same as delivering 10 hours late. In this mindset, an on-time delivery is considered a success while a slightly late delivery is considered a mediocre success and the very late delivery is considered a failure. We account for slightly late deliveries by defining a penalty function that is multiplied by the probabilities of late delivery. The penalty is a function of lateness and starts from 1 corresponding to lateness of zero and diminishes to zero as the lateness increases. We explain this further in section \ref{formulation}.

The formulation is presented here:



\subsection{Operators}

$\times_{k} :$  $\mathscr{k}$-mode tensor product

$\odot :$ Hadamard product

$\oslash :$ Hadamard division

$\bigvee :$ max operator

$\bigwedge :$ min operator


\subsection{Parameters}

$s \in S:$ index for suppliers , $\quad n_S = |S|:$ number of suppliers

$d \in D:$ index for depots , $\quad n_D = |D|:$ number of depots

$f \in F:$ index for storage facilities , $\quad n_F = |F|:$ number of storage facilities

$a \in A:$ index for destinations , $\quad n_A = |A|:$ number of destination nodes


$e \in E:$ index for aid commodity , $\quad n_E = |E|:$ number of aid commodities

$\boldsymbol{\O}_{m,n}:$ $m \times n$ matrix of zeros

$\zeta :$ tuning parameter for the $\psi$ function

$\eta^c :$ probability threshold for preprocessing viable transfers $\in (0,1)$

$\eta^h :$ probability threshold for homotopy algorithm $\in (0,1)$

$\kappa, \theta :$ shape and scale parameters of a Gamma distribution

$\mu, \sigma :$ mean and standard deviation of a distribution (hrs)

$\nu :$ load factor for vehicles

$\psi :$ late delivery penalty function

$b^a:$ amount of aid available at each location (ton)

$b^h:$ shortage of commodities on each level used for optimization (ton)

$b^n:$ amount of demand for aid at the destinations (ton)

$b^u:$ amount of unsatisfied demand at the destinations (ton)


$c^t:$ cost of travel between the nodes ($\$$)

$f:$ probability density function

$\boldsymbol{J}_{m,n}:$ $m \times n$ matrix of ones

$k^c:$ capacity consumed by aid commodities (ton or $m^3$)

$k^v:$ capacity of vehicles (ton or $m^3$)

$m:$ total number of vehicles on a level in the network

$n:$ total number of nodes on a level in the network

$P(X \leq Y):$ probability of random variable $X$ being less than $Y$

$p^b:$ probability estimate used for optimizing binary variables $\in (0,1)$

$p^c:$ probability estimate used for preprocessing $\in (0,1)$

$p^d:$ probability of successful transfer of aid between two levels $\in (0,1)$

$\mathcal{R} :$ reliability measure $\in (0,100)$

$\mathcal{R}^u :$ deficiencies in the reliability measure $\in (0,100)$

$t^a:$ time that each vehicle becomes available to be dispatched (hrs past current time)

$t^{i}:$ available time intervals between two levels (hrs)

$t^{ic}:$ available time intervals used in preprocessing (hrs)

$t^l:$ lateness of delivery compared to $t^n$ (hrs)

$t^n:$ time demand for delivery of aid at destinations (hrs past current time)


$\boldsymbol{U^c}:$ network connectivity matrix between two levels (binary)

$\boldsymbol{U^p}:$ parallel dispatch matrix (binary)

$\boldsymbol{U^v}:$ vehicle-node association matrix (binary)

$w:$ importance factor for types of aid

$z:$ available operational budget ($\$$)

\subsection{Variables}

$b^d:$ amount of aid to be dispatched (ton)

$t^d:$ time of dispatch (hrs past current time)

$\boldsymbol{U^d}:$ vehicle-destination dispatch matrix (binary)

\subsection{Formulae} \label{formulation}

The reliability of the system depends on the successful performance of transfers throughout the network. We first calculate the probabilities of successful transfers in different parts of the network with equations (\ref{eq-p1}) - (\ref{eq-p3}).

$\boldsymbol{P}_{n_E,m_F}^{d,(F\rightarrow A)}$ as calculated by equation (\ref{eq-p1}) is the probability of successful delivery of the aid dispatched from level $F$ to level $A$. The rows of $\boldsymbol{P}_{n_E,m_F}^{d,(F\rightarrow A)}$ correspond to the commodities of aid ($n_E$) and its columns ($m_F$) correspond to the vehicles available on level $F$. \begin{equation} \label{eq-p1}
\begin{split}
\boldsymbol{P}_{n_E,m_F}^{d,(F\rightarrow A)} \big( &\boldsymbol{U}_{m_F,n_A}^{d,(F\rightarrow A)}, \boldsymbol{t}_{1,m_F}^{d,(F \rightarrow A)} \big) \\
&= \mathlarger{P} \big( X \leq \boldsymbol{T}_{n_E,m_F}^{i,(F \rightarrow A)} \big) + \mathlarger{P} \big( X > \boldsymbol{T}_{n_E,m_F}^{i,(F \rightarrow A)} \big) . \psi
\big( X - \boldsymbol{T}_{n_E,m_F}^{i,(F \rightarrow A)} \big) \\
&= \mathlarger{\int}_{0}^{\boldsymbol{T}_{n_E,m_F}^{i,(F \rightarrow A)}} f \Big(t|\Big(\big(\{\boldsymbol{\mu}_{m_F,n_A}^{(F\rightarrow A)}, \boldsymbol{\sigma}_{m_F,n_A}^{(F\rightarrow A)}\}\odot \boldsymbol{U}_{m_F,n_A}^{d,(F\rightarrow A)}\big) \times_{n_A} \boldsymbol{J}_{n_A,n_E} \Big)^T \Big) dt \\
&+ \mathlarger{\int}_{\boldsymbol{T}_{n_E,m_F}^{i,(F \rightarrow A)}}^{+\infty} \psi \big( t - \boldsymbol{T}_{n_E,m_F}^{i,(F \rightarrow A)} \big) \\
& \quad \qquad \qquad \odot f \Big(t|\Big(\big(\{\boldsymbol{\mu}_{m_F,n_A}^{(F\rightarrow A)}, \boldsymbol{\sigma}_{m_F,n_A}^{(F\rightarrow A)}\}\odot \boldsymbol{U}_{m_F,n_A}^{d,(F\rightarrow A)}\big) \times_{n_A} \boldsymbol{J}_{n_A,n_E} \Big)^T \Big) dt \\
\end{split}
\end{equation}

where $\boldsymbol{T}_{n_E,m_F}^{i,(F \rightarrow A)}$ is the time intervals between the dispatch times and the demand time which is calculated using equation (\ref{eq-interval}). 
\begin{equation} \label{eq-interval}
\quad \boldsymbol{T}_{n_E,m_F}^{i,(F \rightarrow A)} = \boldsymbol{T}_{n_E,n_A}^{n} \times_{n_A} \boldsymbol{U}_{m_F,n_A}^{{d,(F\rightarrow A)}^T} - \boldsymbol{J}_{n_E,1} \times_1 \boldsymbol{t}_{1,m_F}^{d,(F \rightarrow A)}
\end{equation}

$\boldsymbol{U}_{m_F,n_A}^{d,(F\rightarrow A)}$ and $\boldsymbol{t}_{1,m_F}^{d,(F \rightarrow A)}$ in equation (\ref{eq-p1}) are both variables that influence the $\boldsymbol{P}_{n_E,m_F}^{d,(F\rightarrow A)}$. $\boldsymbol{U}_{m_F,n_A}^{d,(F\rightarrow A)}$ is the vehicle-destination dispatch matrix describing which vehicles on level $F$ should be dispatched to which destination nodes on level $A$. Its rows correspond to the vehicles ($m_F$) and its columns correspond to the destination nodes ($n_A$). $\boldsymbol{t}_{1,m_F}^{d,(F \rightarrow A)}$ describes the time that the aforementioned dispatches should occur on level $F$. $\boldsymbol{T}_{n_E,n_A}^{n}$ is the time that each aid commodity is needed on level $A$. $\boldsymbol{\mu}_{m_F,n_A}^{(F\rightarrow A)}$ and $\boldsymbol{\sigma}_{m_F,n_A}^{(F\rightarrow A)}$ are the mean and standard deviation of the transfer times for each vehicle on level $F$ to arrive at each of the nodes on level $A$. 

$\boldsymbol{P}_{n_E,m_F}^{d,(F\rightarrow A)}$ is the summation of two integrals. The first integral calculates the probability that aid dispatched from level $F$ gets delivered on level $A$ before the demand deadline for the aid commodities, i.e. on-time. The second integral calculates the penalized probability of late delivery. The penalized probabilities are calculated by convolution of probability density function and a penalty function. The probability density function $f(t)$ is discussed further in section \ref{density} and the late delivery penalty function $\psi(t)$ is explained in section \ref{latepenalty}.

$\boldsymbol{P}_{m_D,m_F}^{d,(D\rightarrow F)}$ as calculated by equation (\ref{eq-p2}) is the probability that dispatched aid from level $D$ arrives in time on level $F$ to successfully make the connecting transfer towards level $A$.
\begin{equation} \label{eq-p2}
\begin{split}
\boldsymbol{P}_{m_D,m_F}^{d,(D \rightarrow F)} \big( \boldsymbol{t}_{1,m_F}^{d,(F \rightarrow A)}, \boldsymbol{t}_{1,m_D}^{d,(D \rightarrow F)} \big) 
 = \mathlarger{P} \big( X \leq \boldsymbol{J}_{m_D,1} \times_1 \boldsymbol{t}_{1,m_F}^{d,(F \rightarrow A)} - \boldsymbol{t}_{1,m_D}^{d,(D \rightarrow F)} \times_1 \boldsymbol{J}_{1,m_F} \big) \\
 = \mathlarger{\int}_{0}^{\boldsymbol{J_{m_D,1}} \times_1 \boldsymbol{t_{1,m_F}^{d,(F \rightarrow A)}} - \boldsymbol{t_{1,m_D}^{d,(D \rightarrow F)}} \times_1 \boldsymbol{J_{1,m_F}}}  f \big( t|\{\boldsymbol{\mu}_{m_D,n_F}^{(D \rightarrow F)}, \boldsymbol{\sigma}_{m_D,n_F}^{(D \rightarrow F)} \} \times_{n_F} \boldsymbol{U}_{m_F,n_F}^{v,(F \rightarrow A)} \big) dt
\end{split}
\end{equation}

$\boldsymbol{P}_{m_S,m_D}^{d,(S\rightarrow D)}$ as calculated by equation (\ref{eq-p3}) is the probability of successful transfer of aid from level $S$ to level $D$ so that it can make the connecting transfer towards level $F$.
\begin{equation} \label{eq-p3}
\begin{split}
\boldsymbol{P}_{m_S,m_D}^{d,(S \rightarrow D)} \big( \boldsymbol{t}_{1,m_D}^{d,(D \rightarrow F)} , \boldsymbol{t}_{1,m_S}^{d,(S \rightarrow D)} \big)
= \mathlarger{P} \big( X \leq \boldsymbol{J}_{m_S,1} \times_1 \boldsymbol{t}_{1,m_D}^{d,(D \rightarrow F)} - \boldsymbol{t}_{1,m_S}^{d,(S \rightarrow D)} \times_1 \boldsymbol{J}_{1,m_D} \big) \\ 
= \boldsymbol{\int}_{0}^{\boldsymbol{J_{m_S,1}} \times_1 \boldsymbol{t_{1,m_D}^{d,(D \rightarrow F)}} - \boldsymbol{t_{1,m_S}^{d,(S \rightarrow D)}} \times_1 \boldsymbol{J_{1,m_D}}}  f \big( t| \{\boldsymbol{\mu}_{m_S,n_D}^{(S \rightarrow D)}, \boldsymbol{\sigma}_{m_S,n_D}^{(S \rightarrow D)} \} \times_{n_D} \boldsymbol{U}_{m_D,n_D}^{v,(S \rightarrow D)} \big) dt
\end{split}
\end{equation}

The vehicle-node association matrices $\boldsymbol{U}^{v}$ describe which vehicles are (or will be) available for dispatch at each of the nodes. Vehicles might not be available immediately, so $\boldsymbol{t}^a$ expresses the time when each vehicle is expected to be available. This gives the model the flexibility to account for the vehicles that are currently in service but will become available again later. $\boldsymbol{U}^{c}$ describes the edges in the network, which implies that all nodes in the network are not necessarily connected. The characteristics of the demand are expressed with $\boldsymbol{b}^{n}$ and $\boldsymbol{t}^{m}$. A single node on level $A$ can have multiple demand requests with different demand characteristics.

The total demand on level $A$ can be obtained by equation (\ref{totaldemand}). The importance factor $\boldsymbol{w_{1,n_E}}$ distinguishes the importance of delivery between different aid commodities and allows us to differentiate between different types of aid regarding their unit weight, e.g. one unit weight of medicine might be of equal importance to ten unit weights of food.
\begin{equation} \label{totaldemand}
\lambda = \boldsymbol{w}_{1,n_E} \times_{n_E} \boldsymbol{B}_{n_E,n_A}^{n} \times_{n_A} \boldsymbol{J}_{n_A,1}
\end{equation}

For ease of computation, we require $\boldsymbol{w_{1,n_E}}$ to be normalized as in equation (\ref{wnorm}).
\begin{equation} \label{wnorm}
\boldsymbol{w}_{1,n_E} \times_{n_E} \boldsymbol{J}_{n_E,1} = 1
\end{equation}

The objective function as calculated by equation (\ref{obj}) maximizes the reliability of deliveries according to the demand and subject to the constraints in equations (\ref{c1})-(\ref{c10}):
\begin{equation} \label{obj}
\begin{split}
\max\limits_{\boldsymbol{U}^d,t^{d},b^d} \mathcal{R} & = \frac{100}{\lambda} \; \Bigg[ \color{Mahogany} \boldsymbol{w}_{1,n_E} \times_{n_E} \big( \boldsymbol{B}_{n_E,m_F}^{d,(F \rightarrow A)} \odot \boldsymbol{P}_{n_E,m_F}^{d,(F\rightarrow A)} \big) \times_{m_F} \boldsymbol{J}_{m_F,1} \\
& + \color{Green} \boldsymbol{w}_{1,n_E} \times_{n_E} \bigg( \Big( \big( \boldsymbol{\mathfrak{B}}_{n_E,m_D,m_F}^{d,(D \rightarrow F)} \odot ( \boldsymbol{P}_{1,m_D,m_F}^{d,(D \rightarrow F)} \times_1 \boldsymbol{J}_{1,n_E}) \big) \; {\times}_{m_D} \; \boldsymbol{J}_{1,m_D} \Big) \\
& \color{Green} \odot \boldsymbol{P}_{n_E,m_F}^{d,(F\rightarrow A)} \bigg) \times_{m_F} \boldsymbol{J}_{m_F,1} \\
& + \color{Blue} \boldsymbol{w}_{1,n_E} \times_{n_E} \Bigg( \bigg( \Big( \boldsymbol{\mathfrak{B}}_{n_E,m_S,m_D,m_F}^{d, (S \rightarrow D)} \odot \big( \boldsymbol{P}_{1,m_S,m_D,1}^{d,(S \rightarrow D)} \times_{1} \boldsymbol{J}_{1,n_E} \times_{1} \boldsymbol{J_{1,m_F}} \big) \\
& \color{Blue} \odot \big( \boldsymbol{P}_{1,1,m_D,m_F}^{d,(D \rightarrow F)} \times_{1} \boldsymbol{J}_{1,n_E} \times_{1} \boldsymbol{J}_{1,m_S} \big) \Big) \times_{m_S} \boldsymbol{J}_{1,m_S} \times_{m_D} \boldsymbol{J}_{1,m_D} \bigg) \\
& \color{Blue} \odot \boldsymbol{P}_{n_E,m_F}^{d,(F \rightarrow A)} \Bigg) \times_{m_F} \boldsymbol{J}_{m_F,1} \color{black} \Bigg]
\end{split}
\end{equation}

The first \textcolor{Mahogany}{red} term in equation (\ref{obj}) accounts for the aid commodities that are pre-positioned at the level $F$ and can be dispatched towards level $A$ on demand. The \textcolor{Mahogany}{$\boldsymbol{w_{1,n_E}}$} is the importance factor for each type of aid commodity. \textcolor{Mahogany}{$\boldsymbol{B}_{n_E,m_F}^{d,(F \rightarrow A)}$} is the weight of commodities being dispatched from level $F$ to level $A$ and similarly \textcolor{Mahogany}{$\boldsymbol{P}_{n_E,m_F}^{(F \rightarrow A)}$} is the probability of on-time delivery of aid commodities at level $A$. The second \textcolor{Green}{green} term in equation (\ref{obj}) accounts for the aid commodities that are present on level $D$ and can be sent first to level $F$ and then to level $A$. Similarly the third \textcolor{Blue}{blue} term accounts for the aid commodities that can be obtained on level $S$ and then shipped through the network to level $D$, then to level $F$ and ultimately to level $A$. All these commodities will flow through the network according to our decision variables and each transfer will have a probability for its successful completion.

The first constraint expressed by equation (\ref{c1}) ensures the dispatched aid commodities does not surpass the demand at the destinations.
\begin{equation} \label{c1}
\begin{gathered}
\big( \boldsymbol{B}_{n_E,m_F}^{d,(F \rightarrow A)} + \boldsymbol{\mathfrak{B}}_{n_E,m_D,m_F}^{d,(D \rightarrow F)} \times_{m_D} \boldsymbol{J}_{1,m_D} + \boldsymbol{\mathfrak{B}}_{n_E,m_S,m_D,m_F}^{d,(S \rightarrow D)} \times_{m_D} \boldsymbol{J}_{1,m_D} \times_{m_S} \boldsymbol{J}_{1,m_S} \big) \\
\times_{m_F} \boldsymbol{U}_{m_F,n_A}^{d,(F \rightarrow A)} \leq \boldsymbol{B}_{n_E,n_A}^n
\end{gathered}
\end{equation}
The formulation makes a distinction between the aid that actually exists on a level and the aid that is expected to arrive to the level. The aid commodities dispatched from a level should not be more than the existing commodities on that level. This is imposed by constraints (\ref{c21}) - (\ref{c23}) corresponding to existing commodities on levels $F$, $D$ and $S$ respectively. 
\begin{equation} \label{c21}
\begin{gathered}
\boldsymbol{B}_{n_E,m_F}^{d,(F \rightarrow A)} \times_{m_F} \boldsymbol{U}_{m_F,n_F}^{v,(F \rightarrow A)} \leq \boldsymbol{B}_{n_E,n_F}^{a,(F \rightarrow A)} \\
\end{gathered}
\end{equation}
\begin{equation} \label{c22}
\begin{gathered}
\big( \boldsymbol{\mathfrak{B}}_{n_E,m_D,m_F}^{d,(D \rightarrow F)} \times_{m_D} \boldsymbol{U}_{m_D,n_D}^{v,(D \rightarrow F)} \big) \times_{m_F} \boldsymbol{J}_{1,m_F} \leq \boldsymbol{B}_{n_E,n_D}^{a,(D \rightarrow F)} \\
\end{gathered}
\end{equation}
\begin{equation} \label{c23}
\begin{gathered}
 \big( \boldsymbol{\mathfrak{B}}_{n_E,m_S,m_D,m_F}^{d,(S \rightarrow D)} \times_{m_S} \boldsymbol{U}_{m_S,n_S}^{v,(S \rightarrow D)} \big) \times_{m_F} \boldsymbol{J}_{1,m_F} \times_{m_D} \boldsymbol{J}_{1,m_D} \leq \boldsymbol{B}_{n_E,n_S}^{a,(S \rightarrow D)}
\end{gathered}
\end{equation}
Constraints (\ref{c31}) - (\ref{c33}) ensure that the dispatched aid matrices $B^{d}$ correspond with the dispatch matrices $U^{d}$. In other words, commodities can only be transferred and delivered via dispatched vehicles.
\begin{equation} \label{c31}
\begin{gathered}
\boldsymbol{B}_{n_E,m_F}^{d,(F \rightarrow A)} \odot \big( \boldsymbol{J}_{n_E,n_A} \times_{n_A} \boldsymbol{U}_{m_F,n_A}^{d,(F \rightarrow A)} \big) \leq \boldsymbol{B}_{n_E,m_F}^{d,(F \rightarrow A)}
\end{gathered}
\end{equation}
\begin{equation} \label{c32}
\begin{gathered}
\big( \boldsymbol{\mathfrak{B}}_{n_E,m_D,m_F}^{d,(D \rightarrow F)} \times_{n_E} \boldsymbol{J}_{n_E,1} \big) \odot \big( \boldsymbol{U}_{m_D,n_F}^{d,(D \rightarrow F)} \times_{n_F} \boldsymbol{U}_{m_F,n_F}^{v,(F \rightarrow A)} \big) \\
\odot \; \big( \boldsymbol{J}_{m_D,n_A} \times_{n_A} \boldsymbol{U}_{m_F,n_A}^{d,(F \rightarrow A)} \big) \leq \boldsymbol{\mathfrak{B}}_{n_E,m_D,m_F}^{d,(D \rightarrow F)} \times_{n_E} \boldsymbol{J}_{n_E,1}
\end{gathered}
\end{equation}
\begin{equation} \label{c33}
\begin{gathered}
\big( \boldsymbol{\mathfrak{B}}_{n_E,m_S,m_D,m_F}^{d,(S \rightarrow D)} \times_{n_E} \boldsymbol{J}_{n_E,1} \big) \odot \big( \boldsymbol{U}_{m_S,n_D,1}^{d,(S \rightarrow D)} \times_{n_D} \boldsymbol{U}_{m_D,n_D}^{v,(D \rightarrow F)} \times_{1} \boldsymbol{J}_{1,m_F} \big) \\
\odot \; \bigg( \Big( \big( \boldsymbol{U}_{1,m_D,n_F}^{d,(D \rightarrow F)} \times_{n_F} \boldsymbol{U}_{m_F,n_F}^{v,(F \rightarrow A)} \big) \odot \big( \boldsymbol{J}_{m_D,n_A} \times_{n_A} \boldsymbol{U}_{m_F,n_A}^{d,(F \rightarrow A)} \big) \Big) \times_{1} \boldsymbol{J}_{1,m_S} \bigg) \\ 
\leq \boldsymbol{\mathfrak{B}}_{n_E,m_S,m_D,m_F}^{d,(S \rightarrow D)} \times_{n_E} \boldsymbol{J}_{n_E,1}
\end{gathered}
\end{equation}
Each vehicle in the system has a capacity for carrying commodities. Constraints (\ref{c41}) - (\ref{c43}) enforce this capacity constraint for vehicles dispatched from each level.
\begin{equation} \label{c41}
\begin{gathered}
\boldsymbol{k}_{1,n_E}^{c} \times_{n_E} \big( \boldsymbol{B}_{n_E,m_F}^{d,(F\rightarrow A)} + \boldsymbol{\mathfrak{B}}_{n_E,m_D,m_F}^{d,(D \rightarrow F)} \times_{m_D} \boldsymbol{J}_{1,m_D} \\ 
+ \boldsymbol{\mathfrak{B}}_{n_E,m_S,m_D,m_F}^{d,(S\rightarrow D)} \times_{m_S} \boldsymbol{J}_{1,m_S} \times_{m_D} \boldsymbol{J}_{1,m_D} \big) \leq \boldsymbol{k}_{1,m_F}^{v,(F\rightarrow A)} \\
\end{gathered}
\end{equation}
\begin{equation} \label{c42}
\begin{gathered}
\boldsymbol{k}_{1,n_E}^{c} \times_{n_E} \big( \boldsymbol{\mathfrak{B}}_{n_E,m_D,m_F}^{d,(D \rightarrow F)} \times_{m_F} \boldsymbol{J}_{1,m_F} \\
+ \boldsymbol{\mathfrak{B}}_{n_E,m_S,m_D,m_F}^{d,(S\rightarrow D)} \times_{m_S} \boldsymbol{J}_{1,m_S} \times_{m_F} \boldsymbol{J}_{1,m_F} \big) \leq \boldsymbol{k}_{1,m_D}^{v,(D\rightarrow F)} \\
\end{gathered}
\end{equation}
\begin{equation} \label{c43}
\begin{gathered}
\boldsymbol{k}_{1,n_E}^{c} \times_{n_E} \boldsymbol{\mathfrak{B}}_{n_E,m_S,m_D,m_F}^{d,(S \rightarrow D)} \times_{m_D} \boldsymbol{J}_{1,m_D} \times_{m_F} \boldsymbol{J}_{1,m_F} \leq \boldsymbol{k}_{1,m_S}^{v,(S\rightarrow D)}
\end{gathered}
\end{equation}
Constraint (\ref{c5}) ensures each vehicle is dispatched at most to one destination.
\begin{equation} \label{c5}
\begin{gathered}
\boldsymbol{U}_{m_F,n_A}^{d,(F \rightarrow A)} \times_{n_A} \boldsymbol{J}_{n_A,1} \leq \boldsymbol{J}_{m_F,1} \quad , \quad \boldsymbol{U}_{m_D,n_F}^{d,(D \rightarrow F)} \times_{n_F} \boldsymbol{J}_{n_F,1} \leq \boldsymbol{J}_{m_D,1} \\ \boldsymbol{U}_{m_S,n_D}^{d,(S \rightarrow D)} \times_{n_D} \boldsymbol{J}_{n_D,1} \leq \boldsymbol{J}_{m_S,1}
\end{gathered}
\end{equation}
Constraint (\ref{c6}) ensures that dispatches comply with the existing edges in the network.
\begin{equation} \label{c6}
\begin{gathered}
\boldsymbol{U}_{m_F,n_A}^{d,(F\rightarrow A)} \leq \boldsymbol{U}_{m_F,n_A}^{c,(F\rightarrow A)} \quad , \quad \boldsymbol{U}_{m_D,n_F}^{d,(D\rightarrow F)} \leq \boldsymbol{U}_{m_D,n_F}^{c,(D\rightarrow F)} \quad , \quad \boldsymbol{U}_{m_S,n_D}^{d,(S\rightarrow D)} \leq \boldsymbol{U}_{m_S,n_D}^{c,(S\rightarrow D)}
\end{gathered}
\end{equation}
Constraint (\ref{c7}) ensures the cost is limited by the available budget.
\begin{equation} \label{c7}
\begin{gathered}
\boldsymbol{J}_{1,m_F} \times_{m_F} \big( \boldsymbol{C}_{m_F,n_A}^{t,(F \rightarrow A)} \odot \boldsymbol{U}_{m_F,n_A}^{d,(F \rightarrow A)} \big) \times_{n_A} \boldsymbol{J}_{n_A,1} \\
+ \; \boldsymbol{J}_{1,m_D} \times_{m_D} \big( \boldsymbol{C}_{m_D,n_F}^{t,(D \rightarrow F)} \odot \boldsymbol{U}_{m_D,n_F}^{d,(D \rightarrow F)} \big) \times_{n_F} \boldsymbol{J}_{n_F,1} \\ 
 + \; \boldsymbol{J}_{1,m_S} \times_{m_S} \big( \boldsymbol{C}_{m_S,n_D}^{t,(S \rightarrow D)} \odot \boldsymbol{U}_{m_S,n_D}^{d,(S \rightarrow D)} \big) \times_{n_D} \boldsymbol{J}_{n_D,1} \leq z
\end{gathered}
\end{equation}
Any vehicle can be scheduled for dispatch only after the time vehicle is expected to be available ($\boldsymbol{t^a}$). This is achieved by constraint (\ref{c8}).
\begin{equation} \label{c8}
\boldsymbol{t}_{1,m_F}^{a,(F\rightarrow A)} \leq \boldsymbol{t}_{1,m_F}^{d,(F\rightarrow A)} \quad , \quad  \boldsymbol{t}_{1,m_D}^{a,(D\rightarrow F)} \leq \boldsymbol{t}_{1,m_D}^{d,(D\rightarrow F)} \quad , \quad \boldsymbol{t}_{1,m_S}^{a,(S\rightarrow D)} \leq \boldsymbol{t}_{1,m_S}^{d,(S\rightarrow D)}
\end{equation}
Constraint (\ref{c9}) ensures all amounts of dispatched commodities are non-negative.
\begin{equation} \label{c9}
0 \leq \boldsymbol{B}_{n_E,m_F}^{d,(F\rightarrow A)} , \boldsymbol{\mathfrak{B}}_{n_E,m_D,m_F}^{d,(D \rightarrow F)} , \boldsymbol{\mathfrak{B}}_{n_E,m_S,m_D,m_F}^{d,(S\rightarrow D)}
\end{equation}
Finally, constraint (\ref{c10}) puts a binary restriction on the dispatch matrices. 
\begin{equation} \label{c10}
\quad \boldsymbol{U}_{m_F,n_A}^{d,(F \rightarrow A)}, \boldsymbol{U}_{m_D,n_F}^{d,(D \rightarrow F)}, \boldsymbol{U}_{m_S,n_D}^{d,(S \rightarrow D)} \in \{0,1\} \quad \text{binary}
\end{equation}








\subsection{Probability distribution of the transfer times} \label{density}

The formulation presented in this paper can incorporate any probability distribution for the travel times as long as it is continuous and differentiable. However, we emphasis that the probability distribution has to be chosen carefully and based on recorded data.

It is obvious that travel time between two distinct points is always positive. Therefore the probability of a transfer to be performed instantly or in negative time should always be computationally zero. It is noteworthy that many probability distributions that are commonly used to model transfer times (such as the normal distribution) assign positive probabilities for negative transfer times. Nevertheless, using such distributions can be computationally accurate, since in practice, the probabilities of negative transfer times are effectively zero to single or double precision. For example the cumulative probability of negative transfer times for a normal distribution that its standard deviation is equal to 10 percent of its mean is $10^{-23}$.

The user of our model has the discretion and responsibility to use appropriate probability distribution for transfer times. Choosing the parameters of probability distribution should be based on recorded data, following proper statistical procedures. Furthermore, some probability distributions such as Gamma, F, Rayleigh and chi-squared distributions, can be defined specifically on the semi-infinite interval of $[0,\infty)$ while having zero probability at 0.



\subsection{Late delivery penalty} \label{latepenalty}
As described before, we recognize a slightly late delivery as mediocre success and penalize its probability to account for the mediocrity of lateness. This requires separate calculation of probabilities of on-time delivery and probabilities of late delivery and therefore, $\boldsymbol{P}_{n_E,m_F}^{d,(F\rightarrow A)}$ calculated by equation (\ref{eq-p1}) is the summation of two integrals.

The late delivery penalty function $\psi$, is generally a nonlinear function of lateness, solely defined for positive values. The value of $\psi$ should be 1 for late delivery of zero (i.e. on-time delivery) and is expected to monotonically decrease with increase in lateness. $\psi$ should eventually reach zero for some lateness that is considered fruitless or unacceptable. 

\begin{figure}[h]
\centering
\includegraphics[scale=0.5]{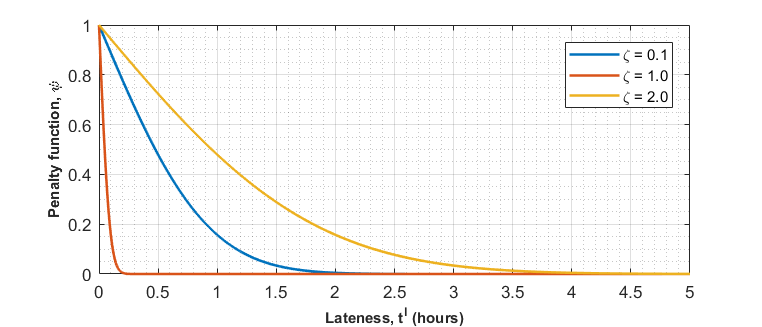}
\vspace{6pt}
\caption{An example of late delivery penalty function $\psi$} 
\label{fig:late}
\end{figure}

Figure~\ref{fig:late} depicts the generic form of the $\psi$ using the error function. In this plot, the penalty function is defined as in equation (\ref{eq-late}) where $t^l$ is the lateness in hours and $\zeta$ is a parameter for tuning the late delivery penalty. Although using the error function is an effective way to define the $\psi$, the user has the discretion to use any other function that is continuous and satisfies the conditions described in the previous paragraph.
\begin{equation} \label{eq-late}
\psi(t^l) = 1-\erf(\frac{t^l}{\zeta})
\end{equation}

As shown in Figure~\ref{fig:late}, decreasing the $\zeta$ will make the penalty more severe and vice verse. Thus, setting the $\zeta$ very small (e.g. 0.001) is equivalent to dropping the second integral in equation (\ref{eq-p1}), implying that late delivery is not valuable even if the lateness is small. We note that because the penalty function is monotonic and strictly smaller than 1, using a more lenient (larger) $\zeta$ does not make the system indifferent with respect to late deliveries. If it is possible to deliver the commodities on-time, the optimal solution will achieve that regardless of the chosen value for $\zeta$.

The trad-off with respect to $\zeta$ becomes meaningful when it is impossible to achieve on-time delivery for all the demand. In such circumstances, choosing a severe penalty will encourage the system not to make dispatches that are expected to be delivered slightly late because those late deliveries will not contribute to the $\mathcal{R}$ and are treated as non-deliveries. Furthermore, depending on the circumstances of a system and the value of $\zeta$, an optimal solution might use a vehicle to deliver a large amount of commodities slightly late instead of using that vehicle to deliver a small amount of commodities on-time.

In conclusion, the $\psi$ function should be defined based on realities on the ground, types of aid commodities and other characteristics of the demand. Thus, it could be plausible to define different $\psi$ functions for each of the nodes on level $A$ and for each of the $n_E$ commodities. In section \ref{algsec}, we morph the $\psi$ function in our homotopic optimization algorithm but the homotopy method should not be confused for the main purpose of the $\psi$ described above.

\subsection{Utilization of the fleet of vehicles}

In this section we look more closely into how the fleet of vehicles are utilized. By using the term \textit{utilization of fleet} we specifically would like to know how the vehicles are dispatched, how much commodities they will carry compared to their capacity and whether there are parallel dispatches in the network.

We first define the term \textit{load factor} of a vehicle as the ratio of the load carried by the vehicle to its capacity. The cost for a vehicle to transfer cargo between two nodes usually does not depend significantly on the load it carries. In other words, when a vehicle is dispatched, we can assume a fixed cost is incurred whether the vehicle is fully loaded or not. Therefore the load factor can be seen as a measure of how efficiently a vehicle is used. Equations (\ref{cu1}) - (\ref{cu3}) calculate the load factor for the vehicles on each level.
\begin{equation} \label{cu1}
\begin{gathered}
\boldsymbol{\nu}_{1,m_F}^{d,(F \rightarrow A)} = \Big( \boldsymbol{k}_{1,n_E}^{c} \times_{n_E} \big( \boldsymbol{B}_{n_E,m_F}^{d,(F\rightarrow A)} + \boldsymbol{\mathfrak{B}}_{n_E,m_D,m_F}^{d,(D \rightarrow F)} \times_{m_D} \boldsymbol{J}_{1,m_D} \\ 
+ \boldsymbol{\mathfrak{B}}_{n_E,m_S,m_D,m_F}^{d,(S\rightarrow D)} \times_{m_S} \boldsymbol{J}_{1,m_S} \times_{m_D} \boldsymbol{J}_{1,m_D} \big)  \Big) \oslash \boldsymbol{k}_{1,m_F}^{v,(F \rightarrow A)}
\end{gathered}
\end{equation}
\begin{equation} \label{cu2}
\begin{gathered}
\boldsymbol{\nu}_{1,m_D}^{d,(D \rightarrow F)} = \Big( \boldsymbol{k}_{1,n_E}^{c} \times_{n_E} \big( \boldsymbol{\mathfrak{B}}_{n_E,m_D,m_F}^{d,(D \rightarrow F)} \times_{m_F} \boldsymbol{J}_{1,m_F} \\
+ \boldsymbol{\mathfrak{B}}_{n_E,m_S,m_D,m_F}^{d,(S\rightarrow D)} \times_{m_S} \boldsymbol{J}_{1,m_S} \times_{m_F} \boldsymbol{J}_{1,m_F} \big) \Big) \oslash \boldsymbol{k}_{1,m_D}^{v,(D \rightarrow F)}
\end{gathered}
\end{equation}
\begin{equation} \label{cu3}
\boldsymbol{\nu}_{1,m_S}^{d,(S \rightarrow D)} = \Big( \boldsymbol{k}_{1,n_E}^{c} \times_{n_E} \boldsymbol{\mathfrak{B}}_{n_E,m_S,m_D,m_F}^{d,(S \rightarrow D)} \times_{m_D} \boldsymbol{J_{1,m_D}} \times_{m_F} \boldsymbol{J}_{1,m_F} \Big) \oslash \boldsymbol{k}_{1,m_S}^{v,(S \rightarrow D)}
\end{equation}

Equation (\ref{cu4}) calculates the average load factor among all dispatched vehicles in the system. The load factor for the non-dispatched vehicles is zero and not included in this average. $\nu^{\text{avg}}$ can be seen as a scalar measure that describes how efficiently vehicles are being used. Because of the constraints (\ref{c41}) - (\ref{c43}), the load factors ($\nu$) will always be upper bounded by 1 which is the logical constraint regarding the capacity of vehicles. $\nu$ is also lower bounded by zero and is continuous between zero and one.

\begin{equation} \label{cu4}
\begin{gathered}
\nu^{\text{avg}} = \bigg( \Big( \boldsymbol{\nu}_{1,m_F}^{d,(F \rightarrow A)} \odot \big( \boldsymbol{J}_{1,n_A} \times_{n_A} \boldsymbol{U}_{m_F,n_A}^{d,(F \rightarrow A)} \big) \Big) \times_{m_F} \boldsymbol{J}_{m_F,1} \\ 
+ \Big( \boldsymbol{\nu}_{1,m_D}^{d,(D \rightarrow F)} \odot \big( \boldsymbol{J}_{1,n_F} \times_{n_F} \boldsymbol{U}_{m_D,n_F}^{d,(D \rightarrow F)} \big) \Big) \times_{m_D} \boldsymbol{J}_{m_D,1} \\ 
+ \Big( \boldsymbol{\nu}_{1,m_S}^{d,(S \rightarrow D)} \odot \big( \boldsymbol{J}_{1,n_D} \times_{n_D} \boldsymbol{U}_{m_S,n_D}^{d,(S \rightarrow D)} \big) \Big) \times_{m_S} \boldsymbol{J}_{m_S,1} \bigg) \\
\oslash \Big( \boldsymbol{U}_{m_F,n_A}^{d,(F \rightarrow A)} \times_{n_A} \boldsymbol{J}_{n_A,1} \times_{m_F} \boldsymbol{J}_{m_F,1} + \boldsymbol{U}_{m_D,n_F}^{d,(D \rightarrow F)} \times_{n_F} \boldsymbol{J}_{n_F,1} \times_{m_D} \boldsymbol{J}_{m_D,1} \\
+ \boldsymbol{U}_{m_S,n_D}^{d,(S \rightarrow D)} \times_{n_D} \boldsymbol{J}_{n_D,1} \times_{m_S} \boldsymbol{J}_{m_S,1} \Big)
\end{gathered}
\end{equation}

We now look into the \textit{parallel dispatches} which refers to two or more vehicles that have the same origin and destination dispatch nodes in the network. $\boldsymbol{U_{m_F,n_A}^{p,(F \rightarrow A)}}$ as calculated by equation (\ref{cu5}) represents the parallel dispatches on level $F$ with entries of 1.
\begin{equation} \label{cu5}
\begin{gathered}
\boldsymbol{U}_{m_F,n_A}^{p,(F \rightarrow A)} = \bigg( \bigvee \Big( \bigwedge \big( \boldsymbol{U}_{m_F,n_A}^{d,(F \rightarrow A)} \times_{m_F} \boldsymbol{U}_{m_F,n_F}^{v,(F \rightarrow A)} , 2 * \boldsymbol{J}_{n_F,n_A} \big) - \boldsymbol{J}_{n_F,n_A} , \boldsymbol{\O}_{n_F,n_A} \Big) \\
\times_{n_F} \boldsymbol{U}_{m_F,n_F}^{v,(F \rightarrow A)} \bigg) \odot \boldsymbol{U}_{m_F,n_A}^{d,(F \rightarrow A)}
\end{gathered}
\end{equation}

Similarly equations (\ref{cu6}) and (\ref{cu7}) identify the parallel dispatches on levels $D$ and $S$, respectively.
\begin{equation} \label{cu6}
\begin{gathered}
\boldsymbol{U}_{m_D,n_F}^{p,(D \rightarrow F)} = \bigg( \bigvee \Big( \bigwedge \big( \boldsymbol{U}_{m_D,n_F}^{d,(D \rightarrow F)} \times_{m_D} \boldsymbol{U}_{m_D,n_D}^{v,(D \rightarrow F)} , 2 * \boldsymbol{J}_{n_D,n_F} \big) - \boldsymbol{J}_{n_D,n_F} , \boldsymbol{\O}_{n_D,n_F} \Big) \\
\times_{n_D} \boldsymbol{U}_{m_D,n_D}^{v,(D \rightarrow F)} \bigg) \odot \boldsymbol{U}_{m_D,n_F}^{d,(D \rightarrow F)}
\end{gathered}
\end{equation}
\begin{equation} \label{cu7}
\begin{gathered}
\boldsymbol{U}_{m_S,n_D}^{p,(S \rightarrow D)} = \bigg( \bigvee \Big( \bigwedge \big( \boldsymbol{U}_{m_S,n_D}^{d,(S \rightarrow D)} \times_{m_S} \boldsymbol{U}_{m_S,n_S}^{v,(S \rightarrow D)} , 2 * \boldsymbol{J}_{n_S,n_D} \big) - \boldsymbol{J}_{n_S,n_D} , \boldsymbol{\O}_{n_S,n_D} \Big) \\
\times_{n_S} \boldsymbol{U}_{m_S,n_S}^{v,(S \rightarrow D)} \bigg) \odot \boldsymbol{U}_{m_S,n_D}^{d,(S \rightarrow D)}
\end{gathered}
\end{equation}

Identifying the parallel dispatches in the system is not necessary in optimizing our model. It merely gives an insight about the inner workings of the system and its optimal solution. For example, it might be useful for pre-disaster planning purposes to know where the parallel dispatches are common in the system.

\section{Optimization Algorithm} \label{algsec}
In this section we propose an optimization algorithm for solving the problem formulated in the previous section. We first investigate the characteristics of formulation and the deficiencies in the reliability and develop a preprocessing procedure to eliminate the enviable edges in the network. Then we decompose the variables and also the network and reason about the optimality of separate components. Finally, we merge the components back together and optimize the system as a whole.

\subsection{Characteristics of the problem}
The formulation is nonlinear, non-smooth and non-concave with mixed-integer variables and convex constraints. Non-concavity requires a global optimization algorithm to be used for solving the problem. Finding a better maximizer can be enormously beneficial since it leads to increasing the reliability of the system and enables the system to deliver more aid on time and with less cost.


The optimization problem in section \ref{formulation} is coded in C++. The linear algebra library named Armadillo (\cite{sanderson2016armadillo}) is used for tensor computations. Although the objective function is non-linear, all the constraints are linear. Furthermore, the objective function is Lipschitz continuous despite being non-smooth. Therefore we can rely on sub-gradients for gradient-based optimization.

Variables are categorizes into three groups: dispatch variables, time variables and commodity variables. These variables are inter-related and optimizing them separately would not lead to an optimal solution for the system. But each category has distinct mathematical characteristics (e.g. binary vs continuous variables) and require different mathematical methods for optimization.

Dispatch variables, which are binary, appear in the objective function and in most of the constraints. The space of feasible solutions for these variables grow rapidly with increase in the number of vehicles and nodes in the network and it would not be possible to investigate all permutations of these variables. Therefore it is necessary to use an integer programming approach to deal with these variables.

Time and commodity variables are continuous and can be optimized using nonlinear programming approaches if binary variables are fixed. But optimizing all the variables together requires a nonlinear mixed-integer optimization algorithm. This category of optimization algorithms is relatively tedious specially for non-concave problems and algorithms are usually designed for specific problems. 

Our main idea for designing an optimization algorithm is based on splitting the binary and continuous variables which is presented in section \ref{secbin}. We then present a homotopy method for optimizing the continuous variables in section \ref{homotopy-section}. Finally in section \ref{mainalg}, we present a comprehensive algorithm for optimizing all the variables in tandem.




\subsection{Preprocessing} \label{sec-preproc}

In pre-processing, we identify and eliminate all the edges in the network that have minute probability of success. The threshold of probability used for this elimination is defined by $\eta^c$. We recommend the default value of $\eta^c$ to be 30 percent but it can be adjusted according to circumstances of the system. $\eta^c$ can also be defined as a tensor with different values for different categories of transfer, but in this formulation, we refer to it as a scalar.

$\boldsymbol{T}_{m_F,n_A}^{ic,(F \rightarrow A)}$ calculated by equation (\ref{eq-tiA}) is the available time interval between the earliest dispatch time $t^a$ for each vehicle on level $F$ and the demand time $t^n$ on level $A$. For multi-commodity systems, the largest demand time will be used in this calculation in order to preserve all edges in the network that have a viable probability of success for at least one commodity.
\begin{equation} \label{eq-tiA}
\boldsymbol{T}_{m_F,n_A}^{ic,(F \rightarrow A)} = \Big( \boldsymbol{J}_{m_F,1} \times_{1} \bigvee\limits_{n_E} \big(\boldsymbol{T}_{n_E,n_A}^{n} \big) - \boldsymbol{t}_{1,m_F}^{{a,(F \rightarrow A)}^T} \times_1 \boldsymbol{J}_{1,n_A} \Big) \odot \boldsymbol{U}_{m_F,n_A}^{c,(F\rightarrow A)}
\end{equation}

Equation (\ref{eq-preprob1}) calculates the $\boldsymbol{P}_{m_F,n_A}^{c,(F\rightarrow A)}$ which represents the highest achievable probability of successful transfers between levels $F$ and $A$.
\begin{equation} \label{eq-preprob1}
\begin{split}
\boldsymbol{P}_{m_F,n_A}^{c,(F\rightarrow A)} &= \mathlarger{P} \big( X \leq \boldsymbol{T}_{m_F,n_A}^{ic,(F \rightarrow A)} \big) + \mathlarger{P} \big( X > \boldsymbol{T}_{m_F,n_A}^{ic,(F \rightarrow A)} \big) . \psi
\big( X - \boldsymbol{T}_{m_F,n_A}^{ic,(F \rightarrow A)} \big) \\
&= \mathlarger{\int}_{0}^{\boldsymbol{T}_{m_F,n_A}^{ic,(F \rightarrow A)}} f \big(t| \{\boldsymbol{\mu}_{m_F,n_A}^{(F\rightarrow A)}, \boldsymbol{\sigma}_{m_F,n_A}^{(F\rightarrow A)} \} \big) dt \\
&+ \mathlarger{\int}_{\boldsymbol{T_{m_F,n_A}^{ic,(F \rightarrow A)}}}^{+\infty} \psi \big( t - \boldsymbol{T}_{m_F,n_A}^{ic,(F \rightarrow A)} \big) \odot f \big(t| \{\boldsymbol{\mu}_{m_F,n_A}^{(F\rightarrow A)}, \boldsymbol{\sigma}_{m_F,n_A}^{(F\rightarrow A)} \} \big) dt \odot \boldsymbol{U}_{m_F,n_A}^{c,(F\rightarrow A)} \\
\end{split}
\end{equation}

Equation (\ref{eq-pre1}) then performs the preprocessing operation on the connectivity matrix between levels $F$ and $A$ of the network. In this equation, $\geq$ is a logical operator that returns a binary tensor.
\begin{equation} \label{eq-pre1}
\boldsymbol{U}_{m_F,n_A}^{c,(F \rightarrow A)} = \boldsymbol{U}_{m_F,n_A}^{c,(F \rightarrow A)} \odot \big( \boldsymbol{P}_{m_F,n_A}^{c,(F \rightarrow A)} \geq \eta^c \big)
\end{equation}

To eliminate enviable edges between levels $D$ and $F$, we have to first calculate the latest time that vehicles on level $F$ can be dispatched so that the minimum probability of $\eta^c$ can be met for at least one destination node. This requires setting $\boldsymbol{P}_{m_F,n_A}^{c,(F \rightarrow A)} = \eta^c$ in equation (\ref{eq-preprob1}) and solving for $\boldsymbol{T}_{m_F,n_A}^{ic,(F \rightarrow A)}$ which will be saved as $\boldsymbol{T}_{m_F,n_A}^{ic,(F \rightarrow A)^{\dagger}}$. 
Using that and similar to the lower level, equation (\ref{eq-tiF}) calculates the time interval available between levels $D$ and $F$.
\begin{equation} \label{eq-tiF}
\begin{gathered}
\boldsymbol{T}_{m_D,n_F}^{ic,(D\rightarrow F)} = \Bigg( \boldsymbol{J}_{m_D,1} \times_{1} \bigvee\limits_{m_F} \Bigg( \bigg( \bigvee\limits_{n_A} \Big( \boldsymbol{J}_{m_F,1} \times_{1} \bigvee\limits_{n_E} \big(\boldsymbol{T}_{n_E,n_A}^{n} \big) - \boldsymbol{T}_{m_F,n_A}^{ic,(F \rightarrow A)^{\dagger}} \Big) \times_{1} \boldsymbol{J}_{1,n_F} \bigg) \\
\odot \boldsymbol{U}_{m_F,n_F}^{v,(F\rightarrow A)} \Bigg) - \boldsymbol{t}_{1,m_D}^{{a,(D \rightarrow F)}^T} \times_1 \boldsymbol{J}_{1,n_F} \Bigg) \odot \boldsymbol{U}_{m_D,n_F}^{c,(D\rightarrow F)}
\end{gathered}
\end{equation}

We now calculate the $\boldsymbol{P}_{m_D,n_F}^{c,(D\rightarrow F)}$ using equation (\ref{eq-preprob2}) which is the highest achievable probability of success for the transfers between levels $D$ and $F$.
\begin{equation} \label{eq-preprob2}
\begin{gathered}
\boldsymbol{P}_{m_D,n_F}^{c,(D \rightarrow F)} = \mathlarger{P} \big( X \leq \boldsymbol{T}_{m_D,n_F}^{ic,(D\rightarrow F)} \big) = \mathlarger{\int}_{0}^{\boldsymbol{T}_{m_D,n_F}^{ic,(D\rightarrow F)}}  f \big( t|\{\boldsymbol{\mu}_{m_D,n_F}^{(D \rightarrow F)}, \boldsymbol{\sigma}_{m_D,n_F}^{(D \rightarrow F)} \} \big) dt \odot \boldsymbol{U}_{m_D,n_F}^{c,(D\rightarrow F)}
\end{gathered}
\end{equation}

And perform the preprocessing between levels $D$ and $F$ using equation (\ref{eq-pre2}).
\begin{equation} \label{eq-pre2}
\boldsymbol{U}_{m_D,n_F}^{c,(D \rightarrow F)} = \boldsymbol{U}_{m_D,n_F}^{c,(D \rightarrow F)} \odot \big( \boldsymbol{P}_{m_D,n_F}^{c,(D \rightarrow F)} \geq \eta^c \big)
\end{equation}

Similar approach can be taken to perform the preprocessing between levels $S$ and $D$ which is left to the reader. These computations do not deal with the variables and are considerably cheaper than each iteration of the optimization algorithm. Therefore, preprocessing can be computationally effective since it eliminates enviable solutions from the space of binary variables.


\subsection{Deficiencies in $\mathcal{R}$ and opportunities for improvement}

We refer to $\mathcal{R}$ as \textit{deficient} when it is less than 100, i.e. imperfect. It follows that any deficiency in $\mathcal{R}$ can be described in terms of the variables and the constraints. Here we categorize the deficiencies into three main groups:

\begin{enumerate}

\item If the total available commodities in the system are less than the total demand, it would lead to deficiency in $\mathcal{R}$. In such circumstances, if there are sufficient time and sufficient fleet of vehicles available to deliver all the commodities in the system on-time with probability 1, the maximum achievable $\mathcal{R}$ can be obtained from equation~\eqref{eq-ceiling}. This is a scalar that is upper-bounded by 100 and defines the ceiling for the global maximizer of the objective function.
\begin{equation} \label{eq-ceiling}
\begin{gathered}
\mathcal{R}^{\text{ceiling}} = 100 \; \boldsymbol{w}_{1,n_E} \times_{n_E} \bigwedge \Big( \big( \boldsymbol{B}_{n_E,n_F}^{a,(F \rightarrow A)} \times_{n_F} \boldsymbol{J}_{n_F,1} + \boldsymbol{B}_{n_E,n_D}^{a,(D \rightarrow F)} \times_{n_D} \boldsymbol{J}_{n_D,1} \\
+ \boldsymbol{B}_{n_E,n_S}^{a,(S \rightarrow D)} \times_{n_S} \boldsymbol{J}_{n_S,1}
\big) \oslash \big( \boldsymbol{B}_{n_E,n_A}^{n} \times_{n_A} \boldsymbol{J}_{n_A,1} \big) , \boldsymbol{J}_{n_E,1} \Big)
\end{gathered}
\end{equation}

\item Deficiency in $\mathcal{R}$ can also come from the dispatch times ($t^d$). If the demand at a destination is urgent and there is not enough time available to deliver the commodities, then the probabilities of on-time delivery would be low and even the deficient $\mathcal{R}$ from previous paragraph would not be achievable. Because of the intertwined and hierarchical relationship among the probabilities in the system, finding the optimal value of dispatch times can be challenging even in a small network.

\item The other adverse influence on $\mathcal{R}$ can arise from the available fleet of vehicles and how they get dispatched in the network. Each vehicle can be dispatched to certain nodes in the network, has a certain capacity to carry commodities and will incur a certain cost depending on where it gets dispatched to. In some nodes, there might be multiple vehicles available to make a certain transfer and there may exist various routes to transfer the commodities from level $S$ to level $A$. If the vehicles are not abundant in the system, it might be challenging or impossible to find a dispatch configuration that can deliver all the commodities with high probability. It might also be impossible to deliver the available commodities because of other constraints such as the available budget. Therefore for a deficient optimal solution, any of the constraints (\ref{c21})-(\ref{c23}), (\ref{c41})-(\ref{c43}) and (\ref{c7}) that is binding, can be described as a cause for deficiency.

\end{enumerate}

Identifying the root causes of deficiencies in $\mathcal{R}$, leads us to find opportunities for improvement. During the optimization process of the system, we are interested to know where the unsatisfied demand are located and what are the possible routes that can possibly accommodate those demand.

Given a solution (a dispatch configuration along with dispatch times and dispatch commodities), unsatisfied demand is basically the slack in constraint (\ref{c1}). Therefore in the ideal case, this constraint should be completely active and binding. Equation (\ref{eq-unsat}) formally defines the slack in constraint (\ref{c1}) as unsatisfied demand $\boldsymbol{B}_{n_E,n_A}^{u}$.
\begin{equation} \label{eq-unsat}
\begin{gathered}
\boldsymbol{B}_{n_E,n_A}^{u} = \boldsymbol{B}_{n_E,n_A}^{n} - \big( \boldsymbol{B}_{n_E,m_F}^{d,(F \rightarrow A)} + \boldsymbol{\mathfrak{B}}_{n_E,m_D,m_F}^{d,(D \rightarrow F)} \times_{m_D} \boldsymbol{J}_{1,m_D} + \boldsymbol{\mathfrak{B}}_{n_E,m_S,m_D,m_F}^{d,(S \rightarrow D)} \\
\times_{m_D} \boldsymbol{J}_{1,m_D} \times_{m_S} \boldsymbol{J}_{1,m_S} \big) \times_{m_F} \boldsymbol{U}_{m_F,n_A}^{d,(F \rightarrow A)}
\end{gathered}
\end{equation}

$\boldsymbol{B_{n_E,n_A}^{u}}$ represents how much of the demand is unsatisfied at each destination node assuming all the dispatched commodities will be successfully delivered at the destinations. But this might not be realistic since it neglects the probabilities. Hence, equation (\ref{eq-defisit}) calculates how the deficiencies in $\mathcal{R}$ are related to the demand at destinations.
\begin{equation} \label{eq-defisit}
\begin{split}
\boldsymbol{\mathcal{R}}_{n_E,n_A}^{u} =\;& 100 \; \big( \boldsymbol{J}_{1,n_A} \times_1 \boldsymbol{w}_{1,n_E} \big) \odot \Bigg[ \Bigg( \bigg( \boldsymbol{B}_{n_E,m_F}^{d,(F \rightarrow A)} + \big( \boldsymbol{\mathfrak{B}}_{n_E,m_D,m_F}^{d,(D \rightarrow F)} \\
&\odot ( \boldsymbol{P}_{1,m_D,m_F}^{d,(D \rightarrow F)} \times_1 \boldsymbol{J_{1,n_E}}) \big) \; \times_{m_D} \; \boldsymbol{J}_{1,m_D} + \Big( \boldsymbol{\mathfrak{B}}_{n_E,m_S,m_D,m_F}^{d, (S \rightarrow D)} \\
&\odot \big( \boldsymbol{P}_{1,m_S,m_D,1}^{d,(S \rightarrow D)} \times_{1} \boldsymbol{J}_{1,n_E} \times_{1} \boldsymbol{J}_{1,m_F} \big) \\
& \odot \; \big( \boldsymbol{P}_{1,1,m_D,m_F}^{d,(D \rightarrow F)} \times_{1} \boldsymbol{J}_{1,n_E} \times_{1} \boldsymbol{J}_{1,m_S} \big) \Big) \times_{m_S} \boldsymbol{J}_{1,m_S} \times_{m_D} \boldsymbol{J}_{1,m_D} \bigg) \\
& \odot \; \boldsymbol{P}_{n_E,m_F}^{d,(F\rightarrow A)} \Bigg) \times_{m_F} \boldsymbol{U}_{m_F,n_A}^{d,(F \rightarrow A)} \Bigg] \oslash \boldsymbol{B}_{n_E,n_A}^{n}
\end{split}
\end{equation}

These parameters will be used later in the optimization algorithm.

\subsection{Dispatch variables} \label{secbin}

Dispatch variables in the formulation are all binary while all other variables are continuous and non-negative. If we only had the binary/dispatch variables, the problem would have been very similar to a knapsack problem where the size of the knapsack is the budget, used up by each vehicle dispatch and the reward is the achieved $\mathcal{R}$. There are yet two caveats to make this idea work:
\begin{enumerate}
\item We have numerous linear constraints and unlike the knapsack problem, the size of knapsack is not the only limitation.
\item for a given solution of binary variables, the reward (achievable $\mathcal{R}$) is not known to us unless we optimize the time and commodity variables.
\end{enumerate}

The first caveat above can be overcome by using a tailored shortest path algorithm that considers a directed acyclic decision graph with edges representing each possible dispatch decision and vertices representing the state of the system. We can define the edges in the decision graph in such a way that constraints (\ref{c5})-(\ref{c6}) cannot be violated. Constraint (\ref{c10}) will also be satisfied in such setting and therefore all constraints related to optimizing the binary variables will be implicitly satisfied except the budget constraint (\ref{c7}) which will be treated as the constraint for the size of knapsack.

Regarding the second caveat, at each vertex on the decision graph, we have a specific dispatch configuration which defines all the binary variables. With the binary variables known at a given vertex, we can optimize the time and commodity variables to obtain the maximum achievable $\mathcal{R}$. In our shortest path mindset, the length of each directed edge will be the negated difference between the maximized $\mathcal{R}$ obtained at each end of the edge. 

By defining the start point in the decision graph as the configuration where no vehicle is being dispatched, we will be searching for the vertex that has the shortest path (i.e. least negated $\mathcal{R}$) that does not violate the budget constraint. This is commonly referred to as \textit{single-source shortest path problem with non-positive weight} in the literature.

We can now build the algorithm \ref{alg-binary} for optimizing the binary variables. Similar to any shortest path algorithm, at each vertex we save the path leading to that vertex (which will be the dispatch configuration), the cost of dispatches (exhausted budget), the maximized $\mathcal{R}$ at the vertex and the corresponding optimal time and commodity variables. When arriving at a new vertex, the first calculation is to check if the budget constraint is satisfied. If not, we should step back to the previous vertex and find a new vertex to explore. But if the budget constraint is satisfied, we will explore the vertex and proceed to optimizing time and commodity variables using a homotopy method as described later in section \ref{homotopy-section}, algorithm \ref{alg-homo}. We refer to this sub-problem optimization as $optimize^h()$, where $h$ stands for homotopy.

Since optimizing the time and commodity variables is expensive, ideally we want to explore the least possible number of vertices in the decision graph and discard enviable paths without actually optimizing the time and commodity variables, rather by drawing inferences from the paths already explored. To achieve this goal, we need to build a strategy based on cheaply computable information and based on previous explorations in the graph.

We first start by drawing the decision graph which is directed and acyclic. The graph will be layered and the number of inner layers will be equal to the available vehicles in the system. The very fist layer will include a single vertex representing the none dispatch decision in the system. This vertex is the starting vertex in the graph and obviously its maximum reliability is $\mathcal{R}=-0$. Considering the direction which commodities flow in the network, we start by the vehicles on level $F$ and then continue to vehicles on levels $D$ and $S$. Each layer in the graph corresponds to one vehicle in the system and has vertices equal to number of distinct destinations that vehicle can be dispatched to plus 1. The plus 1 vertex represents the no dispatch decision for the given vehicle. Therefore the number of edges in the decision graph are equal to $7$ and the number of vertices are $8$. Clearly the number of paths in this decision network will be extremely large and will grow exponentially with number of vehicles. A decision graph will be drawn for one of our numerical examples in section \ref{results}.

\subsubsection{Dispatch decisions on level $F$}

In order to start exploring the decision graph, we first define an easy to solve, binary integer program (BIP) that only considers the vehicles on level $F$. Using equations (\ref{eq-bin1}) and (\ref{eq-bin2}), we calculate the probability estimate $\boldsymbol{P}_{m_F,n_A}^{b,(F\rightarrow A)}$ based on the assumption that all vehicles on level $F$ are dispatched as early as they are available.
\begin{equation} \label{eq-bin1}
\begin{split}
\boldsymbol{T}_{m_F,n_A}^{ib,(F \rightarrow A)} &= \boldsymbol{J}_{m_F,1} \times_{1} \bigg( \Big( \boldsymbol{w_{1,n_E}} \times_{n_E} \big( \boldsymbol{B}_{n_E,n_A}^{n} \odot \boldsymbol{T}_{n_E,n_A}^{n} \big) \Big) \oslash \Big( \boldsymbol{w}_{1,n_E} \times_{n_E} \boldsymbol{B}_{n_E,n_A}^{n} \Big) \bigg) \\
&- \boldsymbol{t}_{1,m_F}^{{a,(F \rightarrow A)}^T} \times_1 \boldsymbol{J}_{1,n_A}
\end{split}
\end{equation}
\begin{equation} \label{eq-bin2}
\begin{split}
\boldsymbol{P}_{m_F,n_A}^{b,(F\rightarrow A)} &= \mathlarger{P} \big( X \leq \boldsymbol{T}_{m_F,n_A}^{ib,(F \rightarrow A)} \big) + \mathlarger{P} \big( X > \boldsymbol{T}_{m_F,n_A}^{ib,(F \rightarrow A)} \big) . \psi
\big( X - \boldsymbol{T}_{m_F,n_A}^{ib,(F \rightarrow A)} \big) \\
&= \mathlarger{\int}_{0}^{\boldsymbol{T}_{m_F,n_A}^{ib,(F \rightarrow A)}} f \big(t| \{\boldsymbol{\mu}_{m_F,n_A}^{(F\rightarrow A)}, \boldsymbol{\sigma}_{m_F,n_A}^{(F\rightarrow A)}\} \big) dt \\
&+ \mathlarger{\int}_{\boldsymbol{T}_{m_F,n_A}^{ib,(F \rightarrow A)}}^{+\infty} \psi \big( t - \boldsymbol{T}_{m_F,n_A}^{ib,(F \rightarrow A)} \big) \odot f \big(t| \{\boldsymbol{\mu}_{m_F,n_A}^{(F\rightarrow A)}, \boldsymbol{\sigma}_{m_F,n_A}^{(F\rightarrow A) }\} \big) dt \odot \boldsymbol{U}_{m_F,n_A}^{c,(F\rightarrow A)}
\end{split}
\end{equation}

We then minimize equation (\ref{bp1-1}) subject to constraints (\ref{bp1-2})-(\ref{bp1-4}) assuming that there is sufficient budget to carry the needed commodities from level $F$ to $A$. This BIP basically finds the dispatch configuration on level $F$ with least cost and highest probability such that adequate capacity is provided to deliver all the demand. Equation (\ref{bp1-1}) divides the cost of each dispatch on level $F$ by the corresponding probability estimate and adds them up. This division by the probability estimate, artificially increases the cost for the dispatches with low chance of success. Hence, vehicles with higher chance of success will be chosen when minimizing the equation (\ref{bp1-1}). If the probability estimate is uniform among all dispatches (e.g. 1), then the BIP reduces to finding the minimum cost dispatch that can provide adequate capacity to satisfy the demand. Constraint (\ref{bp1-2}) ensures adequate capacity is provided for each node on level $A$. Constraint (\ref{bp1-3}) ensures each vehicle is dispatched to one destination and constraint (\ref{bp1-4}) ensures dispatch decisions are binary.
\begin{equation} \label{bp1-1}
\min\limits_{\boldsymbol{U}_{m_F,n_A}^{d,(F \rightarrow A)^{0}}} \boldsymbol{J}_{1,m_F} \times_{m_F} \big( \boldsymbol{C}_{m_F,n_A}^{t,(F \rightarrow A)} \odot \boldsymbol{U}_{m_F,n_A}^{d,(F \rightarrow A)^{0}} \oslash \boldsymbol{P}_{m_F,n_A}^{b,(F \rightarrow A)} \big) \times_{n_A} \boldsymbol{J}_{n_A,1} \\
\end{equation}
\begin{center}
subject to: \\
\end{center}
\begin{equation} \label{bp1-2}
\boldsymbol{k_{1,n_E}^{c}} \times_{n_E} \boldsymbol{B}_{n_E,n_A}^{n} \leq \boldsymbol{k}_{1,m_F}^{v,(F \rightarrow A)} \times_{m_F} \boldsymbol{U}_{m_F,n_A}^{d,(F \rightarrow A)^{0}} \\
\end{equation}
\begin{equation} \label{bp1-3}
\boldsymbol{U}_{m_F,n_A}^{d,(F \rightarrow A)^{0}} \times_{n_A} \boldsymbol{J}_{n_A,1} \leq \boldsymbol{J}_{m_F,1} \\
\end{equation}
\begin{equation} \label{bp1-4}
\boldsymbol{U}_{m_F,n_A}^{d,(F \rightarrow A)^{0}} \in \{ 0,1 \} \qquad \text{binary}
\end{equation}

This BIP might be infeasible if there are not sufficient vehicles available on level $F$ to carry all the demand for commodities. In such case, some portions of the demand will remain unsatisfied because of the inadequate fleet of vehicles on level $F$ and we should find the least infeasible solution that is strictly feasible regarding constraints (\ref{bp1-3}) and (\ref{bp1-4}). Such least infeasible solution will dispatch all the vehicles on level $F$ and satisfy the demand to its achievable maximum. This infeasibility implies deficiency in the optimal $\mathcal{R}$ for the system but not infeasibility of $\boldsymbol{U}_{m_F,n_A}^{d,(F \rightarrow A)^{0}}$ in the main formulation.

Regardless of feasibility, we note that equations (\ref{bp1-1})-(\ref{bp1-4}) do not consider the available budget and therefore after obtaining the optimal solution, we have to recalculate equation (\ref{bp1-1}) with $\boldsymbol{P}_{m_F,n_A}^{b,(F \rightarrow A)} = \boldsymbol{J}_{m_F,n_A}$ and see if it has surpassed $z$. If so, then we have to change our approach and solve another BIP that can deal with the budget.

We now turn our attention to the case where budget is not sufficient to carry all the demand from level $F$ to $A$. For such circumstances, we define an objective function as in equation (\ref{bp2-1}) that minimizes the unsatisfied demand by finding the optimal dispatch decisions on level $F$. This optimization will be constrained by the previously defined constraints (\ref{bp1-3})-(\ref{bp1-4}) in addition to constraint (\ref{bp2-2}) which ensures compliance with the available budget.
\begin{equation} \label{bp2-1}
\begin{split}
\min\limits_{\boldsymbol{U}_{m_F,n_A}^{d,(F \rightarrow A)^{0}}} \bigvee \Big( &\boldsymbol{k}^c_{1,n_E} \times_{n_E} \boldsymbol{B}_{n_E,n_A}^{n} - \boldsymbol{k}_{1,m_F}^{v,(F \rightarrow A)} \\ 
&\times_{m_F} \big( \boldsymbol{U}_{m_F,n_A}^{d,(F \rightarrow A)^{0}} \odot \boldsymbol{P}_{m_F,n_A}^{b,(F \rightarrow A)} \big) , \boldsymbol{\O}_{1,n_A} \Big) \times_{n_A} \boldsymbol{J}_{n_A,1}
\end{split}
\end{equation}
\begin{center}
subject to: constraints (\ref{bp1-3}),(\ref{bp1-4}) and (\ref{bp2-2}) \\
\end{center}
\begin{equation} \label{bp2-2}
\boldsymbol{J}_{1,m_F} \times_{m_F} \big( \boldsymbol{C}_{m_F,n_A}^{t,(F \rightarrow A)} \odot \boldsymbol{U}_{m_F,n_A}^{d,(F \rightarrow A)^{0}} \big) \times_{n_A} \boldsymbol{J}_{n_A,1} \leq z
\end{equation}

Optimizing the equation (\ref{bp2-1}) is harder than optimizing the equation (\ref{bp1-1}). Nevertheless, it is piece-wise linear and convex and therefore can be solved relatively fast with integer programming methods, much cheaper than even calculating gradient of our main objective function.

Using the approach above, we obtain a feasible solution for $\boldsymbol{U}_{m_F,n_A}^{d,(F \rightarrow A)^{0}}$ which is also viable from a limited viewpoint. We also note that, so far we have not accounted for where the commodities are located in the network and we have solely looked at the vehicles on level $F$.

\subsubsection{Dispatch decisions on levels $D$ and $S$}
If available budget is not binding, we can proceed to optimizing the decision variables on level $D$ and then possibly for level $S$. Based on the dispatch decisions obtained in previous section, we first exhaust the available commodities on level $F$ by assigning them to the vehicles that are decided to be dispatched. This assignment can be done by simply maximizing the linear equation (\ref{bp3-1}) subject to constraints (\ref{bp3-2})-(\ref{bp3-4}). Equation (\ref{bp3-1}) which is linear and continuous maximizes the amount of dispatched commodities while constraints (\ref{bp3-2}), (\ref{bp3-3}) and (\ref{bp3-4}) ensure the dispatched commodities do not surpass the available commodities on level $F$, the demand and the capacity of vehicles, respectively.
\begin{equation} \label{bp3-1}
\max\limits_{\boldsymbol{B}_{n_E,m_F}^{d,(F \rightarrow A)^0}} \boldsymbol{w}_{1,n_E} \times_{n_E} \boldsymbol{B}_{n_E,m_F}^{d,(F \rightarrow A)^0} \times_{m_F} \boldsymbol{J}_{m_F,1}
\end{equation}
\begin{center}
subject to: \\
\end{center}
\begin{equation} \label{bp3-2}
\boldsymbol{B}_{n_E,m_F}^{d,(F \rightarrow A)^0} \times_{m_F} \boldsymbol{U}_{m_F,n_A}^{d,(F \rightarrow A)^{0}} \leq \boldsymbol{B}_{n_E,n_A}^{n}
\end{equation}
\begin{equation} \label{bp3-3}
\boldsymbol{B}_{n_E,m_F}^{{d,(F \rightarrow A)^0}} \times_{m_F} \boldsymbol{U}_{m_F,n_F}^{v,(F \rightarrow A)} \leq \boldsymbol{B}_{n_E,n_F}^{a,(F \rightarrow A)}
\end{equation}
\begin{equation} \label{bp3-4}
\boldsymbol{k}_{1,n_E}^{c} \times_{n_E} \boldsymbol{B}_{n_E,m_F}^{d,(F \rightarrow A)^0} \leq \boldsymbol{k}_{1,m_F}^{v,(F \rightarrow A)} \end{equation}

If the optimal solution for equation (\ref{bp3-1}) is not unique, the solution with the most number of nonzero elements shall be chosen. We then calculate the shortage of commodities on level $F$, by maximizing equation (\ref{bp4-1}) subject to constraints (\ref{bp4-2}) and (\ref{bp4-3}) which ensure dispatched commodities do not surpass the demand and vehicle capacities, respectively. 
\begin{equation} \label{bp4-1}
\max\limits_{\boldsymbol{B}_{n_E,m_F}^{h,(F \rightarrow A)^0}} \boldsymbol{w}_{1,n_E} \times_{n_E} \boldsymbol{B}_{n_E,m_F}^{h,(F \rightarrow A)^0} \times_{m_F} \boldsymbol{J}_{m_F,1}
\end{equation}
\begin{center}
subject to: \\
\end{center}
\begin{equation} \label{bp4-2}
\big( \boldsymbol{B}_{n_E,m_F}^{h,(F \rightarrow A)^0} + \boldsymbol{B}_{n_E,m_F}^{d,(F \rightarrow A)^0} \big) \times_{m_F} \boldsymbol{U}_{m_F,n_A}^{d,(F \rightarrow A)^{0}} \leq \boldsymbol{B}_{n_E,n_A}^{n}
\end{equation}
\begin{equation} \label{bp4-3}
\boldsymbol{k}_{1,n_E}^{c} \times_{n_E} \big( \boldsymbol{B}_{n_E,m_F}^{h,(F \rightarrow A)^0} + \boldsymbol{B}_{n_E,m_F}^{d,(F \rightarrow A)^0} \big) \leq \boldsymbol{k}_{1,m_F}^{v,(F \rightarrow A)} \end{equation}

Now we can treat the shortages on level $F$ as demand (similar to demand on level $A$) and find a viable solution for dispatches on level $D$ using the same approach we used for calculating $\boldsymbol{U}_{m_F,n_A}^{d,(F \rightarrow A)^{0}}$. Similarly we shall continue to level $S$. The relevant equations are left to the reader for the sake of brevity.

\subsubsection{Summary}
The approach for finding a viable feasible solution for the dispatch variables is summarized in algorithm \ref{alg-binary}. This naive approach will give us a viable feasible solution on the decision graph that might be close or far from the optimal solution, depending on the circumstances of the system. Later in section \ref{mainalg}, we develop a strategy to modify this feasible solution, searching for the global maximizer of problem.

\IncMargin{1em}
\begin{algorithm}[h] \caption{Finding a viable feasible solution for dispatch and commodity variables} \label{alg-binary}
\SetAlgoLined
\SetKwData{Left}{left}\SetKwData{This}{this}\SetKwData{Up}{up}
\SetKwInOut{Input}{Input}\SetKwInOut{Output}{Output}
\Input{ The aid delivery system with all its parameters }
\Output{ A feasible solution for dispatch and commodity variables $\boldsymbol{U}_{m_S,n_D}^{d,(S \rightarrow D)^0}$, $\boldsymbol{U}_{m_D,n_F}^{d,(D \rightarrow F)^0}$, $\boldsymbol{U}_{m_F,n_A}^{d,(F \rightarrow A)^0}$, $\boldsymbol{B}_{n_E,m_F}^{d,(F \rightarrow A)^0}$ , $\boldsymbol{\mathfrak{B}}_{n_E,m_D,m_F}^{d,(D \rightarrow F)^0}$ and $\boldsymbol{\mathfrak{B}}_{n_E,m_S,m_D,m_F}^{d,(S \rightarrow D)^0}$ }
 Initialize all the dispatch and commodity variables to zeros\;
 Calculate probability estimate $\boldsymbol{P}_{m_F,n_A}^{b,(F\rightarrow A)}$ using equation (\ref{eq-bin1})\; \label{al1PbF}
 Find $\boldsymbol{U}_{m_F,n_A}^{d,(F \rightarrow A)^0}$ by minimizing equation (\ref{bp1-1}) subject to constraints (\ref{bp1-2})-(\ref{bp1-4})\; \label{al1UFb}
 \If{ \text{BIP defined by equations (\ref{bp1-1})-(\ref{bp1-4}) is infeasible} }{
   Find the least infeasible $\boldsymbol{U}_{m_F,n_A}^{d,(F \rightarrow A)^0}$ that is strictly feasible for constraints (\ref{bp1-3})-(\ref{bp1-4})\;
   }
 \If{ $\boldsymbol{U}_{m_F,n_A}^{d,(F \rightarrow A)^0}$ violates the budget constraint }{
   Find $\boldsymbol{U}_{m_F,n_A}^{d,(F \rightarrow A)^0}$ by minimizing equation (\ref{bp2-1}) subject to constraints (\ref{bp1-3}), (\ref{bp1-4}) and (\ref{bp2-2})\;
   return the dispatch and commodity variables\; } \label{al1UFe}
 Calculate $\boldsymbol{B}_{n_E,m_F}^{d,(F \rightarrow A)^0}$ by maximizing equation (\ref{bp3-1}) subject to constraints (\ref{bp3-2})-(\ref{bp3-4})\; \label{al1BFd}
 Calculate the $\boldsymbol{B}_{n_E,m_F}^{h,(F \rightarrow A)^0}$ by maximizing equation (\ref{bp4-1}) subject to constraints (\ref{bp4-2})-(\ref{bp4-3})\; \label{al1BFh}
 \If{$\boldsymbol{B}_{n_E,m_F}^{h,(F \rightarrow A)^0} \ngtr 0$  }{return the dispatch and commodity variables\;}
 Calculate probability estimate $\boldsymbol{P}_{m_D,n_F}^{b,(D \rightarrow F)}$ similar to line \ref{al1PbF}\;
 Find $\boldsymbol{U}_{m_D,n_F}^{d,(D \rightarrow F)^0}$ using similar approach used for $\boldsymbol{U}_{m_F,n_A}^{d,(F \rightarrow A)^0}$ in lines \ref{al1UFb}-\ref{al1UFe} \;
 Calculate $\boldsymbol{\mathfrak{B}}_{n_E,m_D,m_F}^{d,(D \rightarrow F)^0}$ and $\boldsymbol{B}_{n_E,m_D}^{h,(D \rightarrow F)^0}$ similar to lines \ref{al1BFd}-\ref{al1BFh}\;
 \If{$\boldsymbol{B}_{n_E,m_D}^{h,(D \rightarrow F)^0} \ngtr 0$  }{return the dispatch and commodity variables\;}
 Calculate probability estimate $\boldsymbol{P}_{m_S,n_D}^{b,(S \rightarrow D)}$ similar to line \ref{al1PbF}\;
 Find $\boldsymbol{U}_{m_S,n_D}^{d,(S \rightarrow D)^0}$ using similar approach used for $\boldsymbol{U}_{m_F,n_A}^{d,(F \rightarrow A)^0}$ in lines \ref{al1UFb}-\ref{al1UFe} \;
 Calculate $\boldsymbol{\mathfrak{B}}_{n_E,m_S,m_D,m_F}^{d,(S \rightarrow D)^0}$ similar to line \ref{al1BFd}\;
 return the dispatch and commodity variables\; 
\end{algorithm}

\subsection{Homotopy algorithm for optimizing time and commodity variables} \label{homotopy-section}

In this section we develop a homotopy algorithm for optimizing the time and commodity variables in the system. We first exploit the late delivery penalty function $\psi$ to transform the system into a state where all deliveries are possible to be made reliably. In such transformed system, the probabilities of successful transfers and successful deliveries are strictly greater than an adjustable threshold $\eta^h$ and the global optimal solution can be calculated confidently and relatively cheap. We then transform the $\psi$ function back to its original form and trace the path of optimal solutions in an attempt to finish up at the global optimal solution of the original system. This procedure is known as homotopy or continuation method for global optimization.

For any vehicle that gets dispatched from level $S$, the optimal dispatch time will be equal to its corresponding $t^a$. The dispatch time for vehicles that are not scheduled for dispatch will be set to infinity. Therefore, the optimal time variables on level $S$ can be obtained by equation (\ref{eq-ho1}).
\begin{equation} \label{eq-ho1}
\begin{split}
\boldsymbol{t}_{1,m_S}^{{d,(S \rightarrow D)}^*} = \big( \boldsymbol{J}_{1,n_D} \times_{n_D} &\boldsymbol{U}_{m_S,n_D}^{d,(S \rightarrow D)} \big) \odot \boldsymbol{t}_{1,m_S}^{a,(S \rightarrow D)} \\
&+ \big( \boldsymbol{J}_{1,m_S} - \boldsymbol{J}_{1,n_D} \times_{n_D} \boldsymbol{U}_{m_S,n_D}^{d,(S \rightarrow D)} \big) \odot \boldsymbol{\mathlarger{\infty}}_{1,m_S}
\end{split}
\end{equation}

To ensure successful transfers on level $D$, we solve equation (\ref{eq-ho2}) for $\boldsymbol{T}_{m_S,n_D}^{i,(S \rightarrow D)^{\dagger}}$.
\begin{equation} \label{eq-ho2}
\eta^h = \boldsymbol{\int}_{0}^{\boldsymbol{T}_{m_S,n_D}^{i,(S \rightarrow D)^{\dagger}}}  f \big( t| \{\boldsymbol{\mu}_{m_S,n_D}^{(S \rightarrow D)}, \boldsymbol{\sigma}_{m_S,n_D}^{(S \rightarrow D)} \} \big) dt
\end{equation}

Using equation (\ref{eq-ho3}), $\boldsymbol{t}_{1,m_D}^{d,(D \rightarrow F)^{\dagger}}$ can now be calculated which will be later used as the optimal solution of the transformed system and as the starting point for the homotopy algorithm.
\begin{equation} \label{eq-ho3}
\begin{split}
\boldsymbol{t}_{1,m_D}^{d,(D \rightarrow F)^{\dagger}} = \bigvee \bigg( \bigvee\limits_{m_S} \Big( \big( \boldsymbol{T}_{m_S,n_D}^{i,(S \rightarrow D)^{\dagger}} +  \boldsymbol{J}_{1,n_D} &\times_1 \boldsymbol{t}_{1,m_S}^{a,(S \rightarrow D)} \big) \\
&\odot \boldsymbol{U}_{m_S,n_D}^{d,(S \rightarrow D)} \Big) \times_{n_D} \boldsymbol{U}_{m_D,n_D}^{v,(D \rightarrow F)} , \boldsymbol{t}_{1,m_D}^{a,(D \rightarrow F)} \bigg) 
\end{split}
\end{equation}

Similarly for dispatch times on level $F$, we first solve equation (\ref{eq-ho4}) for $\boldsymbol{T}_{m_D,n_F}^{i,(D \rightarrow F)^{\dagger}}$ and then calculate $\boldsymbol{t}_{1,m_F}^{d,(F \rightarrow A)^{\dagger}}$ using equation (\ref{eq-ho5}).
\begin{equation} \label{eq-ho4}
\eta^h = \boldsymbol{\int}_{0}^{\boldsymbol{T_{m_D,n_F}^{i,(D \rightarrow F)^{\dagger}}}}  f \big( t| \{\boldsymbol{\mu}_{m_D,n_F}^{(D \rightarrow F)}, \boldsymbol{\sigma}_{m_D,n_F}^{(D \rightarrow F)} \} \big) dt
\end{equation}
\begin{equation} \label{eq-ho5}
\begin{split}
\boldsymbol{t}_{1,m_F}^{d,(F \rightarrow A)^{\dagger}} = \bigvee \bigg( \bigvee\limits_{m_D} \Big( \big( \boldsymbol{T}_{m_D,n_F}^{i,(D \rightarrow F)^{\dagger}} +  \boldsymbol{J}_{1,n_F} &\times_1 \boldsymbol{t}_{1,m_D}^{a,(D \rightarrow F)} \big) \\
&\odot \boldsymbol{U}_{m_D,n_F}^{d,(D \rightarrow F)} \Big) \times_{n_F} \boldsymbol{U}_{m_F,n_F}^{v,(F \rightarrow A)} , \boldsymbol{t}_{1,m_F}^{a,(F \rightarrow A)} \bigg) 
\end{split}
\end{equation}

$\boldsymbol{t}_{1,m_S}^{d,(S \rightarrow D)^{*}}$, $\boldsymbol{t}_{1,m_D}^{d,(D \rightarrow F)^{\dagger}}$ and $\boldsymbol{t}_{1,m_F}^{d,(F \rightarrow A)^{\dagger}}$ give us the ideal dispatch times that would ensure successful flow of commodities from level $S$ and successful transfers on levels $D$ and $F$. Based on that, we now need to transform the system to ensure the deliveries on level $A$ are successful as well. This transformation is used by adjusting the late delivery penalty function $\psi$.

As it was explained in section \ref{latepenalty}, the main purpose of the $\psi$ function is to realistically model the consequences of late delivery of commodities. Here we exploit that function to artificially allow for late deliveries and transform the system to an ideal state where all the deliveries can be made reliably, i.e. with probability $\eta^h$. In this transformation, we assume the late delivery penalty is defined by equation (\ref{eq-late}) for each demand via a given tensor $\boldsymbol{\zeta}_{n_E,m_F}$. 

To obtain the ideal system, we need to calculate the required $\zeta$ parameters that would ensure successful delivery on level $A$. Hence, we first calculate $\boldsymbol{T}_{n_E,m_F}^{i,(F \rightarrow A)^{\dagger}}$ by inserting $\boldsymbol{t}_{1,m_F}^{d,(F \rightarrow A)^{\dagger}}$ in equation (\ref{eq-interval}). Then we solve the integral in equation (\ref{eq-ho6}) for $\boldsymbol{\zeta}_{n_E,m_F}^{0}$.
\begin{equation} \label{eq-ho6}
\begin{split}
& \eta^h = \mathlarger{P} \big( X \leq \boldsymbol{T}_{n_E,m_F}^{i,(F \rightarrow A)^{\dagger}} \big) + \mathlarger{P} \big( X > \boldsymbol{T}_{n_E,m_F}^{i,(F \rightarrow A)^{\dagger}} \big) . \psi
\big( X - \boldsymbol{T}_{n_E,m_F}^{i,(F \rightarrow A)^{\dagger}} \big) \\
&= \mathlarger{\int}_{0}^{\boldsymbol{T}_{n_E,m_F}^{i,(F \rightarrow A)^{\dagger}}} f \Big(t|\Big(\big(\{\boldsymbol{\mu}_{m_F,n_A}^{(F\rightarrow A)}, \boldsymbol{\sigma}_{m_F,n_A}^{(F\rightarrow A)}\}\odot \boldsymbol{U}_{m_F,n_A}^{d,(F\rightarrow A)}\big) \times_{n_A} \boldsymbol{J}_{n_A,n_E} \Big)^T \Big) dt \\
&+ \mathlarger{\int}_{\boldsymbol{T}_{n_E,m_F}^{i,(F \rightarrow A)^{\dagger}}}^{+\infty} \bigg( \boldsymbol{J}_{n_E,m_F} - \erf \Big( \big( t - \boldsymbol{T}_{n_E,m_F}^{i,(F \rightarrow A)^{\dagger}}  \big) \oslash \boldsymbol{\zeta}_{n_E,m_F}^{0} \Big) \bigg) \\
& \quad \qquad \qquad \odot f \Big(t|\Big(\big(\{\boldsymbol{\mu}_{m_F,n_A}^{(F\rightarrow A)}, \boldsymbol{\sigma}_{m_F,n_A}^{(F\rightarrow A)}\}\odot \boldsymbol{U}_{m_F,n_A}^{d,(F\rightarrow A)}\big) \times_{n_A} \boldsymbol{J}_{n_A,n_E} \Big)^T \Big) dt
\end{split}
\end{equation}

We point out that equation (\ref{eq-ho6}) might not have a solution if $\eta^h$ is chosen smaller than the first integral on the right hand side. In such cases setting $\eta^h$ very close (or equal) to 1 will guarantee the existence of solution.

To find the starting point for the commodity variables, we can consider the time and dispatch variables as constants and optimize the linear program (LP) in equation (\ref{eq-ho7}) subject to our previously defined constraints (\ref{c1}) - (\ref{c43}) and (\ref{c9}). The space of feasible solutions for this LP will be strictly non-empty (since it contains the all zero solution) and the maximizer obtained from it will define the ceiling for the result of homotopy algorithm.
\begin{equation} \label{eq-ho7}
\begin{split}
\max\limits_{\boldsymbol{b^d}} \mathcal{R} = \frac{1}{\lambda} \; \Big[ \boldsymbol{w}_{1,n_E} \times_{n_E} &\bigg( \boldsymbol{B}_{n_E,m_F}^{d,(F \rightarrow A)} \times_{m_F} \boldsymbol{J}_{m_F,1} \\
& + \eta^h \odot \boldsymbol{\mathfrak{B}}_{n_E,m_D,m_F}^{d,(D \rightarrow F)} {\times}_{m_D} \boldsymbol{J}_{m_D,1} \times_{m_F} \boldsymbol{J}_{m_F,1} \\
& + {\eta^h}^2 \odot \boldsymbol{\mathfrak{B}}_{n_E,m_S,m_D,m_F}^{d, (S \rightarrow D)} \times_{m_S} \boldsymbol{J}_{m_S,1} \times_{m_D} \boldsymbol{J}_{m_D,1} \times_{m_F} \boldsymbol{J}_{m_F,1} \Bigg]
\end{split}
\end{equation}

As with any linear program, the optimal solution obtained from solving equation (\ref{eq-ho7}) might not be unique and/or there might exist degenerate constraints. In such cases, the optimal solution with the least number of non-zero elements should be chosen.

Now, we have all the ingredients to use the homotopy method presented in algorithm \ref{alg-homo}. As mentioned earlier, this homotopy algorithm considers the binary dispatch variables to be constant and does not optimize them. $\boldsymbol{t}_{1,m_S}^{d,(S \rightarrow D)^{*}}$ is also optimal and therefore fixed in this setting.
\IncMargin{1em}
\begin{algorithm} \caption{Homotopy algorithm for optimizing time and commodity variables} \label{alg-homo}
\SetAlgoLined
\SetKwData{Left}{left}\SetKwData{This}{this}\SetKwData{Up}{up}
\SetKwInOut{Input}{Input}\SetKwInOut{Output}{Output}
\Input{ Instance of the aid delivery system with all its parameters and fixed dispatch variables $\boldsymbol{U}_{m_S,n_D}^{d,(S \rightarrow D)}, \boldsymbol{U}_{m_D,n_F}^{d,(D \rightarrow F)} , \boldsymbol{U}_{m_F,n_A}^{d,(F \rightarrow A)}$}
\Output{Optimal values for $\boldsymbol{t}_{1,m_D}^{d,(D \rightarrow F)}$, $\boldsymbol{t}_{1,m_F}^{d,(F \rightarrow A)}$ , $\boldsymbol{B}_{n_E,m_F}^{d,(F \rightarrow A)}$ , $\boldsymbol{\mathfrak{B}}_{n_E,m_D,m_F}^{d,(D \rightarrow F)}$ and $\boldsymbol{\mathfrak{B}}_{n_E,m_S,m_D,m_F}^{d,(S \rightarrow D)}$  }
 Calculate the optimal value for $\boldsymbol{t}_{1,m_S}^{d,(S \rightarrow D)^{*}}$ using equation (\ref{eq-ho1})\;
 Define $optimize^{s}()$ as the function that optimizes equation (\ref{obj}) with fixed dispatch variables, fixed $\boldsymbol{t}_{1,m_S}^{d,(S \rightarrow D)^{*}}$, subject to constraints (\ref{c1})-(\ref{c10}) and returns optimal time and commodity variables\;
 Calculate the starting point for time variables: $\boldsymbol{t}_{1,m_D}^{d,(D \rightarrow F)^{0}}$ and $\boldsymbol{t}_{1,m_F}^{d,(F \rightarrow A)^{0}}$ using equations (\ref{eq-ho2})-(\ref{eq-ho5})\;
 Calculate the starting point for commodity variables: $\boldsymbol{B}_{n_E,m_F}^{d,(F \rightarrow A)^0}$ , $\boldsymbol{\mathfrak{B}}_{n_E,m_D,m_F}^{d,(D \rightarrow F)^0}$ and $\boldsymbol{\mathfrak{B}}_{n_E,m_S,m_D,m_F}^{d,(S \rightarrow D)^0}$ by optimizing the LP in equation (\ref{eq-ho7}) subject to constraints (\ref{c1})-(\ref{c43}) and (\ref{c9})\;
 Save the starting points as $opt^0$\;
 Calculate the transformation parameter $\boldsymbol{\zeta}_{n_E,m_F}^{0}$ by solving equation (\ref{eq-ho6})\;
 \eIf{$\boldsymbol{\zeta}_{n_E,m_F}^{0} \leq \boldsymbol{\zeta}_{n_E,m_F}$}{
 $optimize^s()$ with $\boldsymbol{\zeta}_{n_E,m_F}$ and with $opt^0$ as starting point\;
 \text{return the optimal solution}\;
 }{
  \text{Choose number of iterations} $\iota^h$\;
  \For{$i=1$ \KwTo $\iota^h$}{ 
    $\boldsymbol{\zeta}_{n_E,m_F}^{i} = \boldsymbol{\zeta}_{n_E,m_F}^{0} + \frac{i}{\iota^h} ( \boldsymbol{\zeta}_{n_E,m_F}-\boldsymbol{\zeta}_{n_E,m_F}^{0} )$\; \label{alg-homo-step}
    $optimize^s()$ with $\boldsymbol{\zeta}_{n_E,m_F}^{i}$ and with $opt^{i-1}$ as starting point\;
    save the obtained optimal solution as $opt^i$\;
  }
  \text{return the optimal solution}\;
 }
\end{algorithm}

If the $\boldsymbol{\zeta}_{n_E,m_F}^{0}$ obtained from equation (\ref{eq-ho7}) is smaller than the original $\boldsymbol{\zeta_{n_E,m_F}}$, then the system has ample time to make all the deliveries and therefore the homotopy transformation will be unnecessary. In such case we will optimize the system for all the time and commodity variables using the starting points obtained above. This will count as only one iteration of the homotopy algorithm (the if statement in algorithm \ref{alg-homo}) and the resulting optimal solution will be in close proximity of the starting point.

The $optimize^s()$ function in algorithm \ref{alg-homo} can be based on any gradient-based nonlinear optimization method that considers linear constraints. We strongly recommend trust-region methods based on our numerical experiments. In our examples in section \ref{results}, we have used the NLopt solver package (\cite{nlopt}) for $optimize^s()$.

\subsection{Main algorithm for optimizing the whole problem} \label{mainalg}

We put together algorithms \ref{alg-binary} and \ref{alg-homo} and all the components in previous sections to develop the overall optimization algorithm \ref{alg-main}. In simple words, we use the viable solution we obtained in section \ref{secbin} and try to improve it by modifying the dispatch variables while exploring the decision graph and optimizing the continuous variables at each iteration using homotopy. Modifying the dispatch decisions means changing the destinations of already dispatched variables and add or remove the vehicles from the list of dispatched vehicles.

Here we develop the strategy for this process. At each iteration of algorithm (i.e. after optimizing time and commodities for a given dispatch configuration), we are interested to know the answers to these intertwined questions:
\begin{enumerate}
\item How much each of the vehicles is contributing to the $\mathcal{R}$?
\item For vehicles with negligible contribution, how much of their capacity is used? Is the node they are serving have unsatisfied demand? If a vehicle from upper level brings more commodities for them, will their contribution change?
\item How much of the demand is unsatisfied at each node, as calculated by equation (\ref{eq-unsat})?
\item For each node with unsatisfied demand, how much unused capacity exists on the vehicles already destined to the node?
\item For nodes with unsatisfied demand, which nodes on upper levels have the commodities and the vehicles to satisfy them? 
\end{enumerate}

Based on the answers to these questions, we build our strategy to explore the decision graph and find the shortest path. At each iteration of the algorithm, we shall first alter the dispatch variables and then use homotopy to optimize the continuous variables to find the maximum achievable $\mathcal{R}$ for that dispatch configuration. We will also keep track of the paths explored in the decision graph similar to a shortest path algorithm.

Based on the contribution to $\mathcal{R}$ we can identify the dispatches with negligible contribution. We evaluate the cause for their negligible contributions and perform a local search to decide whether their contribution can be improved. Vehicles that their contribution to $\mathcal{R}$ remain below a threshold $\eta^c$, will be set to no dispatch.

We then look at the unsatisfied demand and the most viable options that can satisfy them. This process will only add vehicles to the dispatch matrix. By repeating these two processes, the dispatch matrix will be altered until we converge to an optimal solution that cannot be improved anymore.

closest vehicle on level F available to be dispatched towards the node - probability - commodity

\begin{equation} \label{eq-sort}
\begin{split}
\boldsymbol{\mathcal{R}}_{n_E,n_A}^{u} &= 100 \; \big( \boldsymbol{J}_{1,n_A} \times_1 \boldsymbol{w}_{1,n_E} \big) \odot \Bigg[ \Bigg( \bigg( \boldsymbol{B}_{n_E,m_F}^{d,(F \rightarrow A)} + \big( \boldsymbol{\mathfrak{B}}_{n_E,m_D,m_F}^{d,(D \rightarrow F)} \\
&\odot ( \boldsymbol{P}_{1,m_D,m_F}^{d,(D \rightarrow F)} \times_1 \boldsymbol{J}_{1,n_E}) \big) \; \times_{m_D} \; \boldsymbol{J}_{1,m_D} + \Big( \boldsymbol{\mathfrak{B}}_{n_E,m_S,m_D,m_F}^{d, (S \rightarrow D)} \\
&\odot \big( \boldsymbol{P}_{1,m_S,m_D,1}^{d,(S \rightarrow D)} \times_{1} \boldsymbol{J}_{1,n_E} \times_{1} \boldsymbol{J}_{1,m_F} \big) \\
& \odot \; \big( \boldsymbol{P}_{1,1,m_D,m_F}^{d,(D \rightarrow F)} \times_{1} \boldsymbol{J}_{1,n_E} \times_{1} \boldsymbol{J}_{1,m_S} \big) \Big) \times_{m_S} \boldsymbol{J}_{1,m_S} \times_{m_D} \boldsymbol{J}_{1,m_D} \bigg) \\
& \odot \; \boldsymbol{P}_{n_E,m_F}^{d,(F\rightarrow A)} \Bigg) \times_{m_F} \boldsymbol{U}_{m_F,n_A}^{d,(F \rightarrow A)} \Bigg] \oslash \boldsymbol{B}_{n_E,n_A}^{n}
\end{split}
\end{equation}

\IncMargin{1em}
\begin{algorithm}[h] \caption{Main optimization algorithm for optimizing all variables} \label{alg-main}
\SetAlgoLined
\SetKwData{Left}{left}\SetKwData{This}{this}\SetKwData{Up}{up}
\SetKwInOut{Input}{Input}\SetKwInOut{Output}{Output}
\Input{ The aid delivery system with all its parameters }
\Output{ Optimal values for all variables }
 Calculate $\boldsymbol{U}_{m_S,n_D}^{d,(S \rightarrow D)^0}$, $\boldsymbol{U}_{m_D,n_F}^{d,(D \rightarrow F)^0}$, $\boldsymbol{U}_{m_F,n_A}^{d,(F \rightarrow A)^0}$, $\boldsymbol{B}_{n_E,m_F}^{d,(F \rightarrow A)^0}$ , $\boldsymbol{\mathfrak{B}}_{n_E,m_D,m_F}^{d,(D \rightarrow F)^0}$ and $\boldsymbol{\mathfrak{B}}_{n_E,m_S,m_D,m_F}^{d,(S \rightarrow D)^0}$ using algorithm \ref{alg-binary}\;
 Define $optimize^{h}()$ as optimizing time and commodity variables for a fixed dispatch configuration using algorithm \ref{alg-homo}\;
 $preprocess()$ the network as described in section \ref{sec-preproc}\;
 list of explored paths\;
 list of vertices in the decision graph and \;
 Initialize iteration counter $i=0$\;
 Initialize \textit{not-optimal = true}\;
 \While{not-optimal}{
 $optimize^{h}()$ with $\boldsymbol{U}_{m_S,n_D}^{d,(S \rightarrow D)^{i}}$, $\boldsymbol{U}_{m_D,n_F}^{d,(D \rightarrow F)^{i}}$ and $\boldsymbol{U}_{m_F,n_A}^{d,(F \rightarrow A)^{i}}$ and save as $opt^i$\;
 contribution\;
 \;
 unsatisfied\;
 \;
 \eIf{$\big( \boldsymbol{U}_{m_S,n_D}^{d,(S \rightarrow D)^{i+1}} = \boldsymbol{U}_{m_S,n_D}^{d,(S \rightarrow D)^{i}} \big) \; \& \; \big( \boldsymbol{U}_{m_D,n_F}^{d,(D \rightarrow F)^{i+1}} = \boldsymbol{U}_{m_D,n_F}^{d,(D \rightarrow F)^{i}} \big) \; \& \; \big( \boldsymbol{U}_{m_F,n_A}^{d,(F \rightarrow A)^{i+1}} = \boldsymbol{U}_{m_F,n_A}^{d,(F \rightarrow A)^{i}} \big)$}{not-optimal = false\;}
 {$i = i + 1$\;}
 }
 return the $opt^i$\;
\end{algorithm}

\subsection{Budget tightening for cost minimization} \label{tightening}

For systems with limited fleet and/or budget resources, the optimization will automatically and implicitly avoid inefficiency in the system. For example, if the cargo assigned to two parallel dispatches can be merged and sent via one of the vehicles with similar probability of success, the optimization module, in search of a better maximizer for $\mathcal{R}$, will naturally consider this reassignment in order to use one of the vehicles (or its assumed cost) for other activities in the system. 

On the other hand, if two parallel dispatches have different probabilities of successful transfer, the optimization module will assign as much commodity as possible to the vehicle with higher probability of success. It follows that any commodities assigned to the vehicle with lower probability of success, are either in excess of the capacity of the other vehicle or the commodities are not available at the time the vehicle with higher probability is dispatched. Therefore, merging the two dispatches is either impossible because of capacity constraints or sub-optimal because it would yield lower reliability.

But when the resources are abundant in the system, the maximum reliability might be obtained with various different dispatch decisions corresponding to different costs. Going back to our two vehicle example, the maximal reliability, might be obtainable whether two dispatches are merged or not when probabilities of success are maximal for both dispatches and the budget constraint is not binding. In such cases, the model will not be able to distinguish between those optimal solutions because they all deliver the same reliability and satisfy all the constraints. Nevertheless the optimal solution we are interested to find is the one that maximizes the $\mathcal{R}$ with minimum cost. Here we explain a bound tightening method for achieving this goal.

Because the secondary objective is to minimize the cost, the best approach would be to gauge the sensitivity of the maximized $\mathcal{R}$ to the budget constraint and try to possibly achieve the same reliability with reduced budget. It is obvious that tightening the budget bound in equation (\ref{c1}) will tighten the space of feasible solutions. Therefore, as the available budget ($z$) is reduced, the optimal $\mathcal{R}$ will either decrease or remain constant.

The tightening method starts with optimizing the system with the original available budget. After finding the initial optimal solution, the available budget will be reduced through a series of iterations. At each iteration the problem will be optimized using the solution obtained at the previous iteration as the starting point. The tightening can continue until $\mathcal{R}$ starts to decrease or when it falls below a certain value.



Another issue that is implicitly addressed with the tightening is related to the end tail of probability distributions. For probability distributions that their integral reach the 1 asymptotically, achieving the exact one depends on the precision of arithmetic. As a result, in our model, when the normalized reliability gets relatively close to 1, boosting it towards 1 might require extensive amount of resources and the formulation will try to achieve that if the available budget and the available fleet allows. This might not be wise and cost effective in many situations and hence it would be wise to always gauge the sensitivity of the obtained reliability with respect to the available budget and possibly consider tightening the budget to achieve a reliability minimally smaller than the one found as a maximizer. This approach is similar to putting a cap on the value of $\mathcal{R}$ by imposing a constraint. For a perfect system that achieves 1, the cap on $\mathcal{R}$ can be 0.99 as demonstrated for our numerical examples.


\section{Results} \label{results}

In order to evaluate the model, here we investigate two examples and present the results.

\subsection{Case 1} \label{case1}

\begin{figure}[h]
\centering
\includegraphics[scale=0.9]{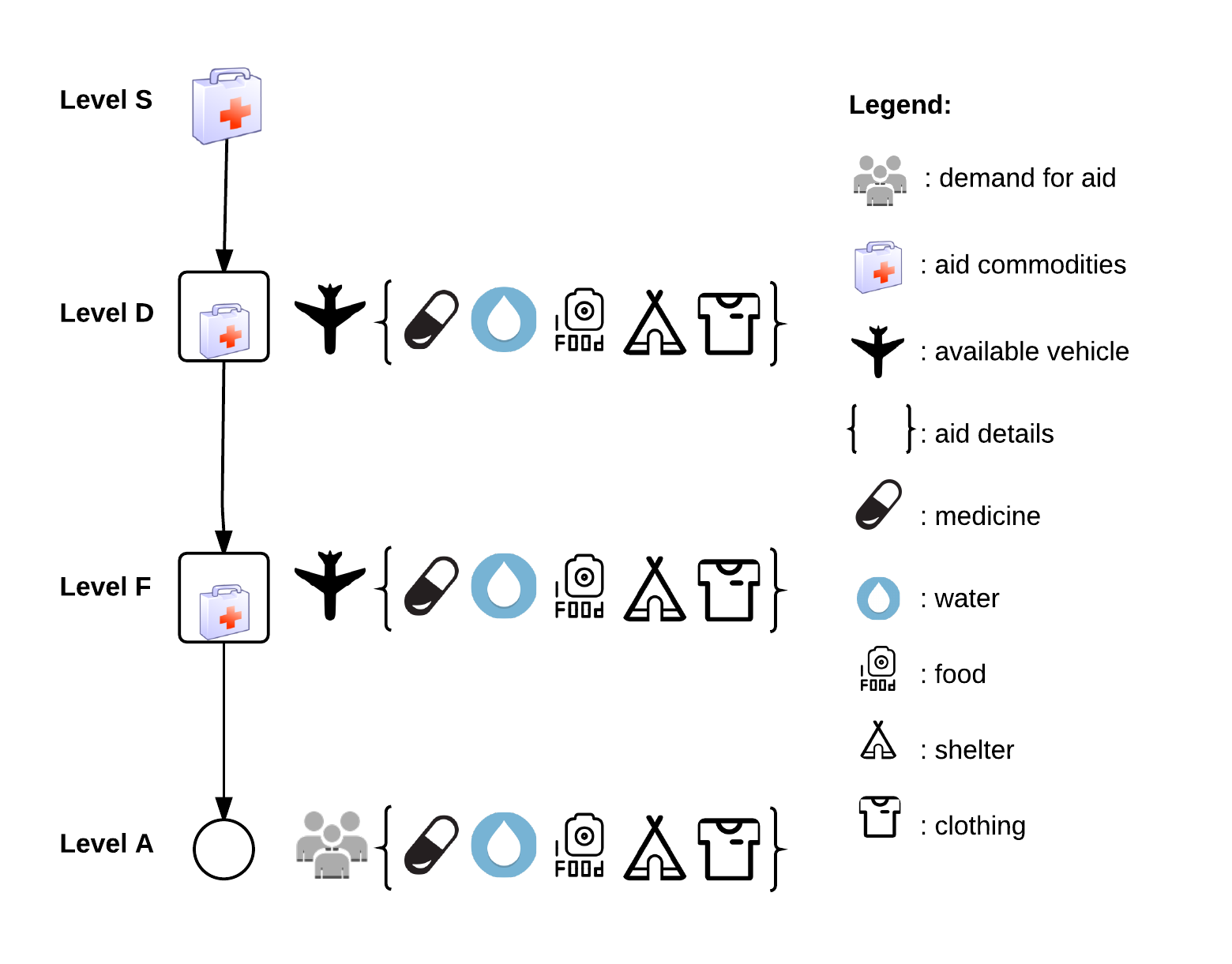}
\caption{Network for Case 1} 
\label{fig:case1}
\end{figure}

The smallest possible network has one node on each level, as shown in Figure~\ref{fig:case1} and we investigate it to gain some insight about the model. We assume that there are only three vehicles available in the system, one on each level ready to be dispatched. There are six variables: $\boldsymbol{B}_{n_E,1}^{d,(F \rightarrow A)}$ and ${t_{1,1}^{d,(F \rightarrow A)}}$ regarding the dispatch from level $F$ to $A$, $\boldsymbol{\mathfrak{B}}_{n_E,1,1}^{d,(D \rightarrow F)}$ and ${t_{1,1}^{d,(D \rightarrow F)}}$ regarding the dispatch from level $D$ to $F$ and finally $\boldsymbol{\mathfrak{B}}_{n_E,1,1}^{d,(S \rightarrow D)}$ and ${t_{1,1}^{d,(S \rightarrow D)}}$ regarding the dispatch from level $S$ to $D$.

We further limit the variables by assuming that the demand consists of five aid commodities and the available commodities on levels $F$, $D$ and $S$ are exactly equal to 40, 20 and 40 percent of the demand, respectively. Moreover, we assume all vehicles have ample capacity to carry all the demand and the available budget is not binding. With these assumptions, the optimal amount of dispatched aid becomes fixed and equal to the available aid on each level. For the sake of demonstration, we consider the transfer time parameters to be as follows according to normal and gamma distributions:

$\mu_{1,1}^{(S \rightarrow D)} = 1.5 \text{ hrs} \quad , \quad \sigma_{1,1}^{(S \rightarrow D)} = 0.3 \text{ hrs} \quad , \quad \mu_{1,1}^{(D \rightarrow F)} = 2.5 \text{ hrs} , \quad , \quad \sigma_{1,1}^{(D \rightarrow F)} = 0.8 \text{ hrs} \quad $

$\kappa_{1,1}^{(F \rightarrow A)} = 4.0 \quad , \quad \theta_{1,1}^{(F \rightarrow A)} = 0.3 \quad , \quad t_{(1:5),1}^{m} = 6.2 \text{ hrs} \quad , \quad \zeta=0.001$

\begin{figure}[htb]
\centering
\begin{tikzpicture}
  \node (img1)  {\includegraphics[scale=0.6]{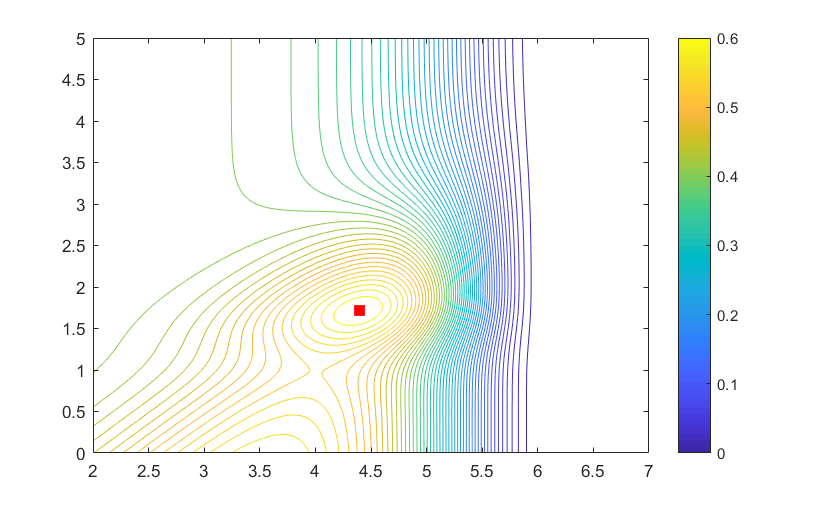}};
  \node[below=of img1, node distance=0cm, yshift=1.6cm,font=\color{black}] {$t_{1,1}^{d,(F \rightarrow A)}$ (hrs)};
  \node[left=of img1, node distance=0cm, rotate=90, anchor=center,yshift=-1.9cm,font=\color{black}] {$t_{1,1}^{d,(D \rightarrow F)}$ (hrs)};
  \end{tikzpicture}
  \vspace{0pt}
  \caption{Contour plot of $\mathcal{R}$ in case 1 with respect to two variables} 
  \label{fig:case1contour}
\end{figure}

$\kappa$ and $\theta$ are shape and scale parameters of gamma distribution used for the transfer time between levels $F$ and $A$. Considering the available time frame, it is viable to transfer the commodities from level $S$ and $D$ to satisfy the demand. Therefore the optimal dispatch time on level $S$ will be zero representing the current time. Now the only variables in the system are two scalars: ${t_{1,1}^{d,(F \rightarrow A)}}$ and ${t_{1,1}^{d,(D \rightarrow F)}}$. Figure~\ref{fig:case1contour} shows the contour plot of the objective function ($\mathcal{R}$) with respect to these two time variables.

It is evident from Figure~\ref{fig:case1contour} that the objective function is non-concave even in this simplified example and a global optimization algorithm is required for optimization. The optimal value of $\mathcal{R}$ in this example (marked by a \textcolor{Red}{red square} on the contour plot in Figure~\ref{fig:case1contour}) is 59.32 and the optimal value for the variables are: $t_{1,1}^{d,(F \rightarrow A)} = 4.399 \text{ hrs}$ and $t_{1,1}^{d,(D \rightarrow F)} = 1.712 \text{ hrs}$.

We now apply the homotopy algorithm \ref{alg-homo} to demonstrate how the homotopy transforms the system and finds the optimal solution. Since the optimal value of commodity variables are known, the homotopy will only perform on the time variables. We choose $\eta^h=1-10^2=0.99$ and obtain $t_{1,1}^{d,(F \rightarrow A)^0} = 4.399$ , $t_{1,1}^{d,(D \rightarrow F)^0} = 1.712$ and $\zeta^0 \approx 70$(rounded up).

\begin{figure}[p]
\begin{subfigure}{0.33\textwidth}
\includegraphics[width=1\linewidth, height=4cm]{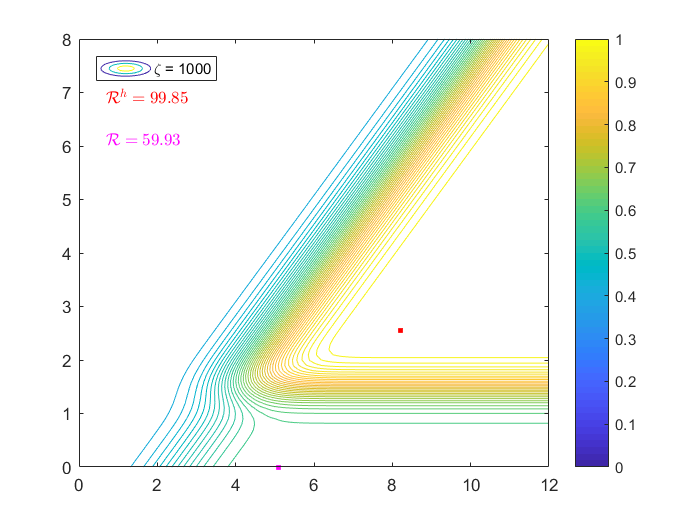}
\caption{}
\label{fig:subim1}
\end{subfigure}
\begin{subfigure}{0.33\textwidth}
\includegraphics[width=1\linewidth, height=4cm]{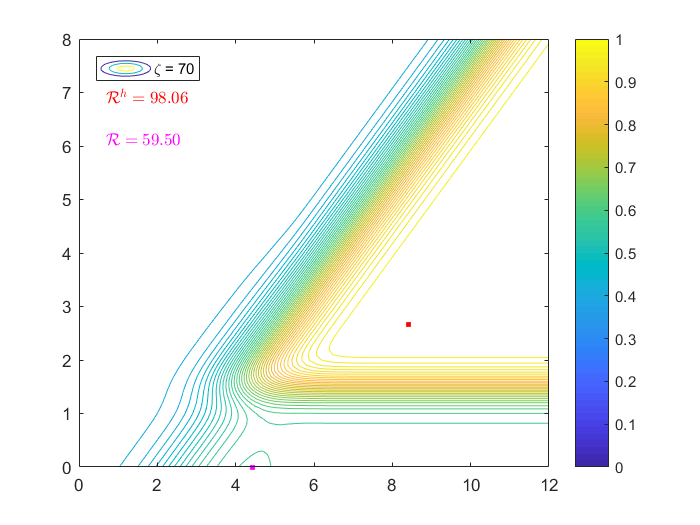}
\caption{}
\label{fig:subim2}
\end{subfigure}
\begin{subfigure}{0.33\textwidth}
\includegraphics[width=1\linewidth, height=4cm]{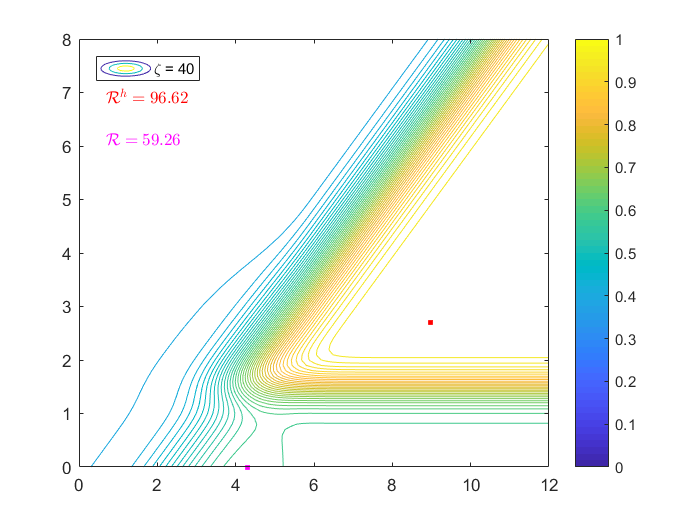}
\caption{}
\label{fig:subim3}
\end{subfigure}
\begin{subfigure}{0.33\textwidth}
\includegraphics[width=1\linewidth, height=4cm]{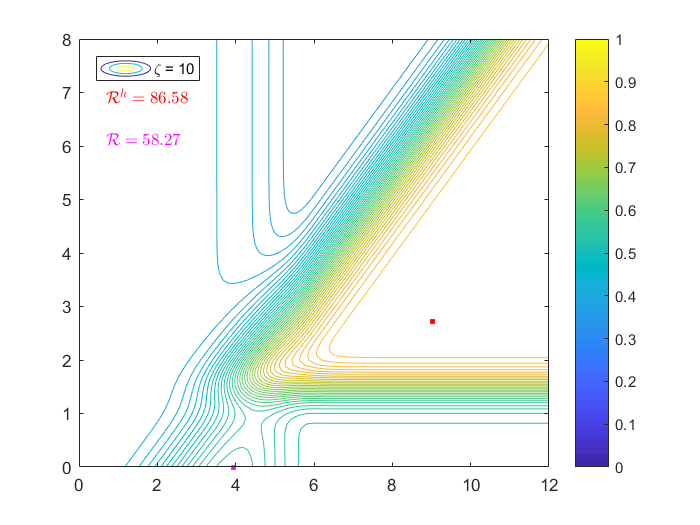} 
\caption{}
\label{fig:subim4}
\end{subfigure}
\begin{subfigure}{0.33\textwidth}
\includegraphics[width=1\linewidth, height=4cm]{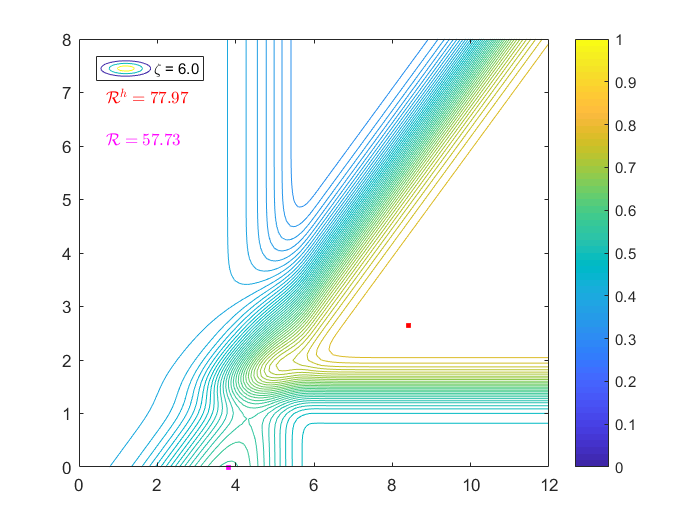}
\caption{}
\label{fig:subim5}
\end{subfigure}
\begin{subfigure}{0.33\textwidth}
\includegraphics[width=1\linewidth, height=4cm]{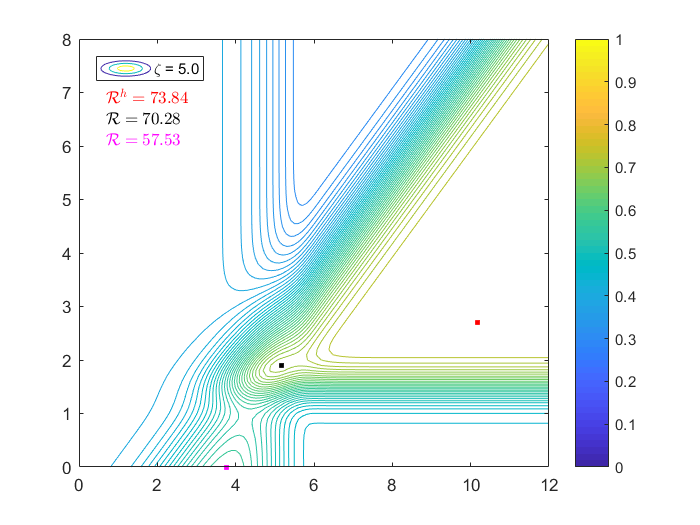}
\caption{}
\label{fig:subim6}
\end{subfigure}
\begin{subfigure}{0.33\textwidth}
\includegraphics[width=1\linewidth, height=4cm]{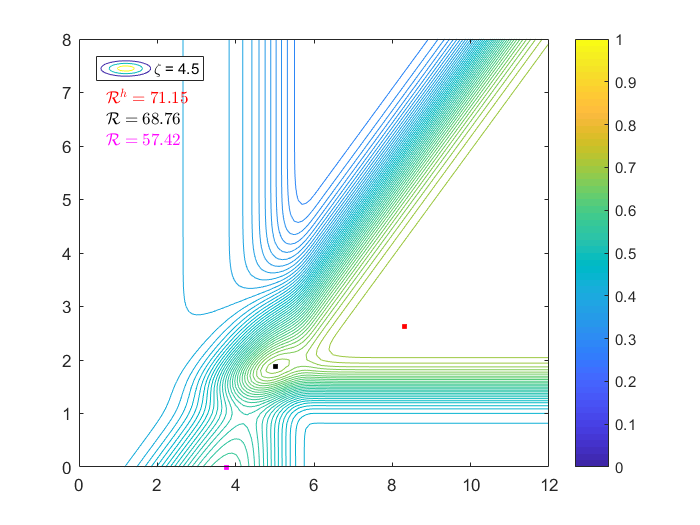}
\caption{}
\label{fig:subim7}
\end{subfigure}
\begin{subfigure}{0.33\textwidth}
\includegraphics[width=1\linewidth, height=4cm]{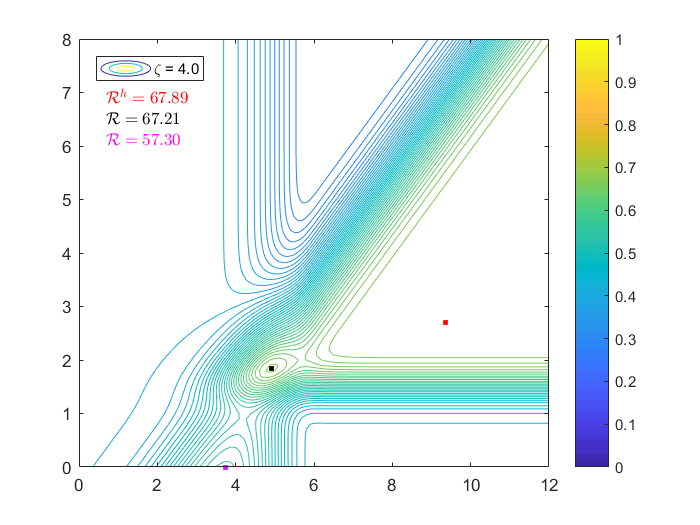}
\caption{}
\label{fig:subim8}
\end{subfigure}
\begin{subfigure}{0.33\textwidth}
\includegraphics[width=1\linewidth, height=4cm]{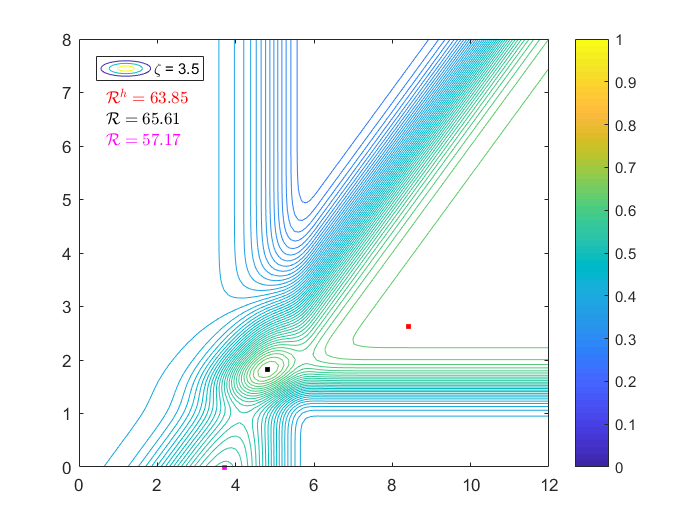}
\caption{}
\label{fig:subim9}
\end{subfigure}
\begin{subfigure}{0.33\textwidth}
\includegraphics[width=1\linewidth, height=4cm]{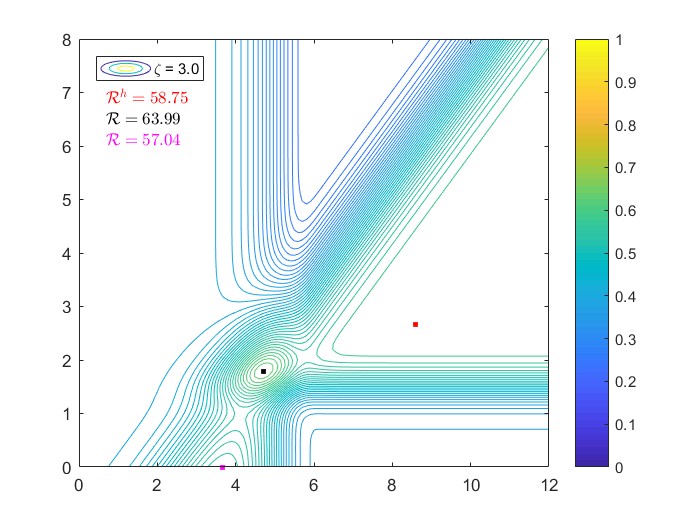}
\caption{}
\label{fig:subim10}
\end{subfigure}
\begin{subfigure}{0.33\textwidth}
\includegraphics[width=1\linewidth, height=4cm]{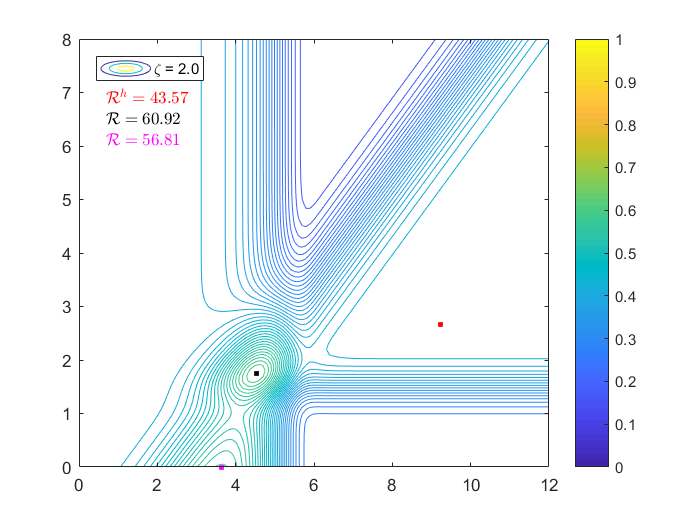}
\caption{}
\label{fig:subim11}
\end{subfigure}
\begin{subfigure}{0.33\textwidth}
\includegraphics[width=1\linewidth, height=4cm]{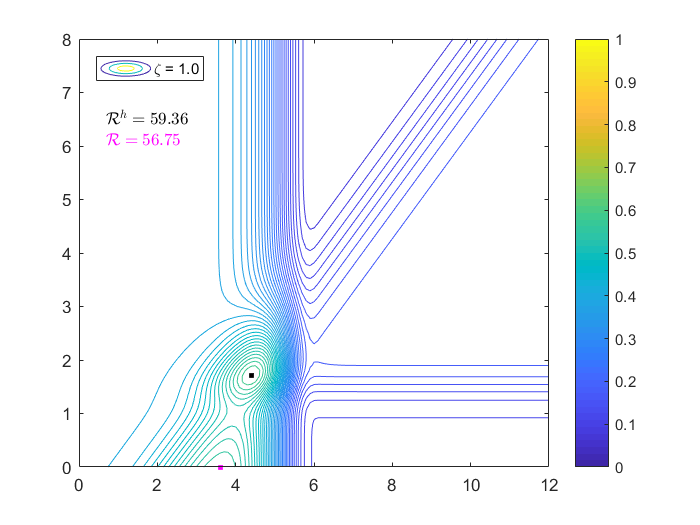}
\caption{}
\label{fig:subim12}
\end{subfigure}
\begin{subfigure}{0.33\textwidth}
\includegraphics[width=1\linewidth, height=4cm]{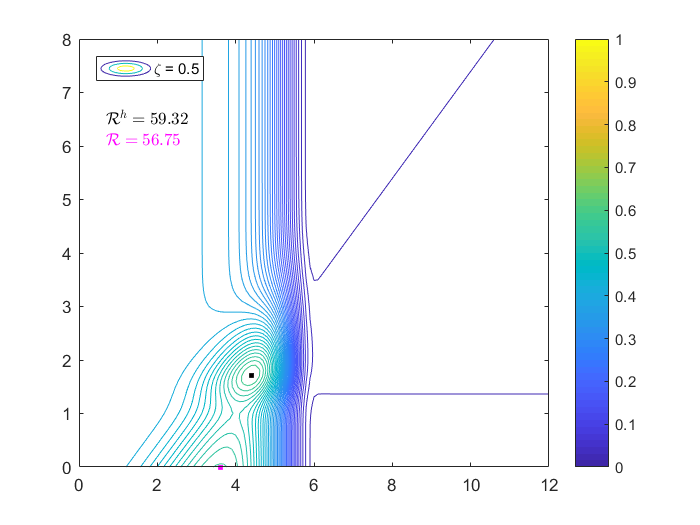}
\caption{}
\label{fig:subim13}
\end{subfigure}
\begin{subfigure}{0.33\textwidth}
\includegraphics[width=1\linewidth, height=4cm]{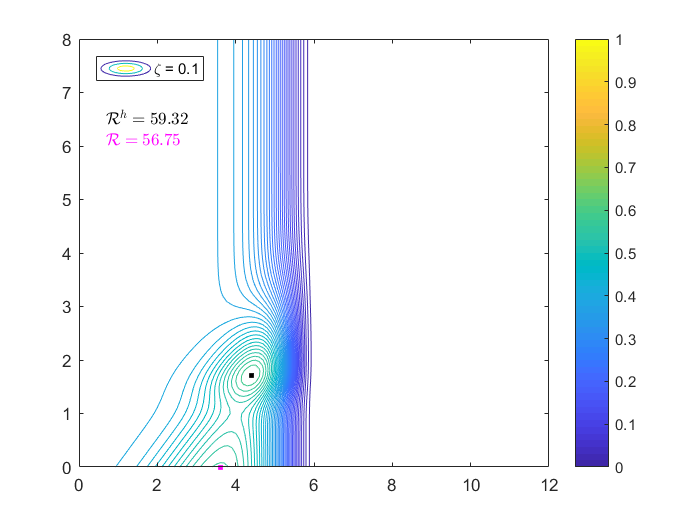}
\caption{}
\label{fig:subim14}
\end{subfigure}
\begin{subfigure}{0.33\textwidth}
\includegraphics[width=1\linewidth, height=4cm]{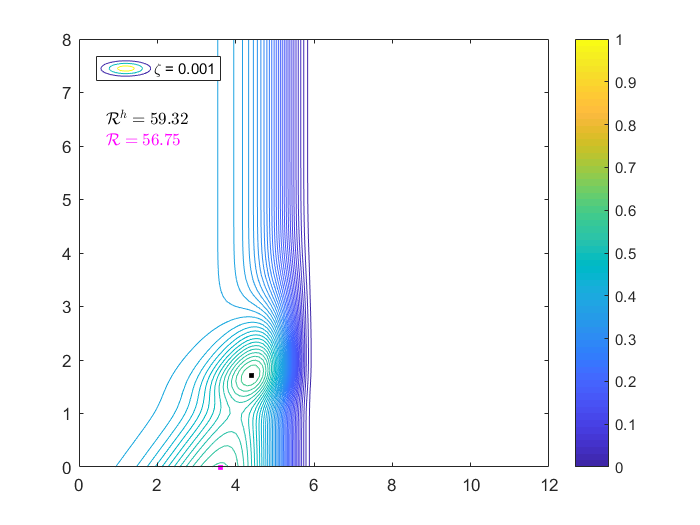}
\caption{}
\label{fig:subim15}
\end{subfigure}
\vspace{8pt}
\caption{Transformation of system in Case 1 using homotopy algorithm \ref{alg-homo}}
\label{fig:case1homotopy}
\end{figure}

We observe $\zeta^0 > \zeta$ and choose the number of iterations $\iota^h = 12$. For the sake of demonstration, instead of using uniform step sizes for $\zeta$ as suggested in line \ref{alg-homo-step} of algorithm \ref{alg-homo}, we choose nonuniform step sizes as shown in Figure~\ref{fig:case1homotopy}.

Figure~\ref{fig:case1homotopy} depicts the contour plot of objective function at each state of transformation in the homotopy algorithm and the \textcolor{red}{red square} shows the location of optimal solution found at each iteration. At each iteration, the optimal solution from the previous iteration is used as the starting point. The algorithm will actually start at subgraph \ref{fig:subim2} and proceed forward, while subgraph \ref{fig:subim1} is presented merely as additional information, showing how the contour plot changes for a very large $\zeta$. The area without a contour in vicinity of the optimal solution, towards the right end of the subgraphs \ref{fig:subim1}-\ref{fig:subim11}, represents a domain where the objective function is locally (and sometimes globally) maximal and constant (flat) in single precision. Therefore, thinking in floating point arithmetic, the flat domain can be seen as a set of optimal solutions and the \textcolor{red}{red square} in each of those subgraphs represents the point from that set with the least euclidean norm.

In subgraphs \ref{fig:subim1}-\ref{fig:subim5}, corresponding to $\zeta$ shrinking from 1000 to 6, there are only two local maximizers. One of these maximizers is the \textcolor{red}{red square} described in previous paragraph and the second one is the single point marked with a \textcolor{magenta}{magenta square}. In subgraph \ref{fig:subim2}, the starting point of the homotopy algorithm coincides with the \textcolor{red}{red square} and is the global maximizer by construction. As we continue decreasing the $\zeta$, at around $\zeta=5$, a new local maximizer emerges which is shown with a \textbf{black square} in subgraph \ref{fig:subim6}. All the three maximizers persist in the contour plots until the vicinity of $\zeta = 1$ corresponding to subgraph \ref{fig:subim12}. At this $\zeta$, the two areas with negative curvature, that correspond to local maximizers marked with \textcolor{red}{red} and \textbf{black} squares, become connected by a continuous negative curvature and therefore the local maximizer with \textcolor{red}{red} square vanishes. The two other maximizers remain near the same location, as we decrease the $\zeta$ to its real value 0.001 in subgraph \ref{fig:subim15}, which is equivalent to Figure~\ref{fig:case1contour}.

It goes without saying that finding the global maximizer of the system in Figure~\ref{fig:case1contour} is not a hard task and most global optimization algorithms can find it even faster than the homotopy method. As we will show in the next case, the effectiveness of homotopy algorithm is observed when solving large networks with numerous variables.

\subsection{Case 2} \label{case2}

\begin{figure}[h]
\includegraphics[scale=.8]{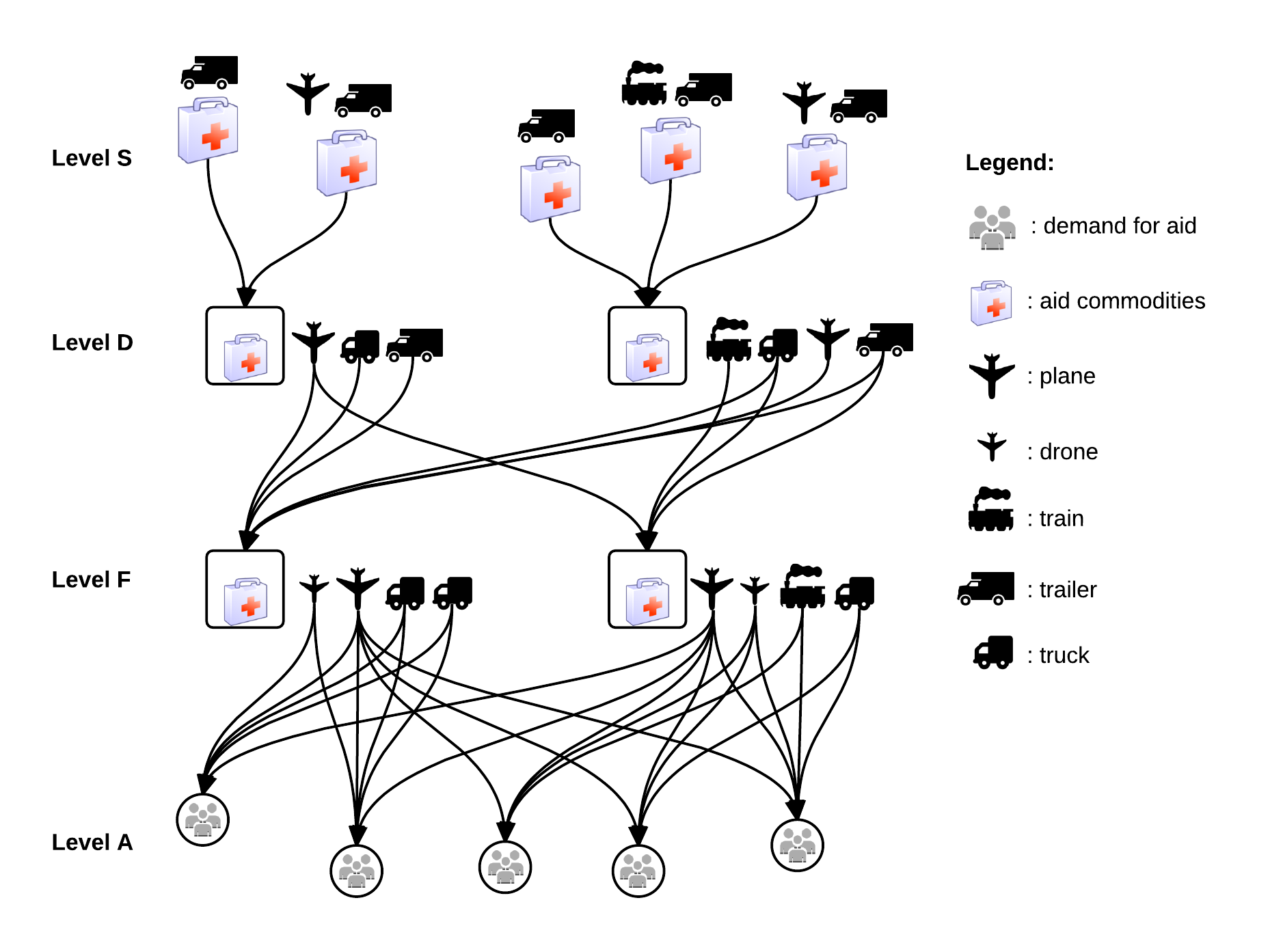}
\centering
\caption{Network for examples of case 2} 
\label{fig:case2}
\end{figure}

In this case we consider a single commodity system as in Figure~\ref{fig:case2} and investigate how different variables and constraints influence the reliability of the system, $\mathcal{R}$.

We first investigate the deficiencies in the system. If the total available commodities in the system are less than the total demand, it would lead to deficiency in $\mathcal{R}$. In this example, for simplicity, we consider the demand for the commodities to be 20 units at each node on level $A$. Since level $A$ has 5 nodes, the total demand will be 100 units. We further assume total amount of commodities available on levels $F$, $D$ and $A$ are 85 units. It follows that if there are sufficient time and sufficient fleet of vehicles available to deliver all the commodities in the system on-time with probability 1, the maximum achievable $\mathcal{R}$ will be equal to 0.85 according to equation (\ref{eq-ceiling}). This defines the ceiling for the global maximizer of the objective function. By artificially setting all the probabilities to 1 and optimizing only over $B^d$ and $U^d$, the optimization module successfully finds this global maximizer.

The other adverse influence on $\mathcal{R}$ can arise from the available fleet of vehicles and how they get dispatched in the network. In this example, we assume the fleet has sufficient capacities to deliver all the commodities and we further assume ample budget is available to achieve the $\mathcal{R}=0.85$.

Another deficiency in $\mathcal{R}$ can be due to the dispatch and demand times ($t^d$ and $t^n$). If the demand at a destination is urgent and there is not enough time available to deliver the commodities, then the probabilities of on-time delivery would be low and even the $\mathcal{R}=0.85$ would not be achievable. As it was explained earlier, finding the optimal value of dispatch times can be challenging even in a small network as in Figure~\ref{fig:case2}.

In this example we generate random $t^n$ among the nodes on level $A$ with mean of 6.5 and standard deviation of 1.0 hours from the current time. We further generate random transfer times (mean and standard deviation) between the nodes (very small compared to the demand time). When the transfer times are very small compared to the demand time, it is possible to find the optimal $t^d$ such that all probabilities of successful transfer are close to 1. In such circumstance, the ceiling for the global maximizer of $\mathcal{R}$ will still be 0.85 which can be found by optimizing over all variables including the $t^d$.

But when the transfer times are not very small, the trade-offs regarding the probabilities will considerably affect the outcome of the system. We now generate random numbers with mean of 1.5 hours for $\mu^{(S \rightarrow D)}$ and for $\mu^{(D \rightarrow F)}$ and $\mu^{(F \rightarrow A)}$ with mean of 2.5 and 1.0 hours, respectively. The standard deviation of these randomly generated numbers are 0.2 hours. We set the $\sigma^{(S \rightarrow D)}$, $\sigma^{(D \rightarrow F)}$ and $\sigma^{(F \rightarrow A)}$ to 10 percent of their corresponding $\mu$.

As a result it will not be possible anymore to make all the probabilities of successful transfers close to 1. When we optimize the system in its new state, the maximum $\mathcal{R}$ is found to be 0.64 which is considerably less than the 0.85. Because of the relatively large number of variables, we no longer can verify if the $\mathcal{R} = 0.64$ is actually the global maximizer of the system but we know it satisfies the optimality conditions of a local maximizer. Also since the optimization module were able to find the global maximizer in previous setting, we can be hopeful that this maximizer might be the global maximizer.

\section{Conclusion} \label{conclusion}

This work, in progress, proposes a framework for optimizing real-time decisions in a hierarchical aid delivery system. Our formulation considers probabilistic transfer times in the system and maximizes the reliability of deliveries with respect to both the amount and time of demand.


\begin{enumerate}
\item The presented framework is computationally efficient and can be solved for real-time purposes. Our formulation is a mixed integer non-linear non-concave maximization problem.

\item We consider a realistic hierarchical network and consider the completion time of activities to be probabilistic. It was shown through numerical examples, that a humanitarian disaster relief system can become robust and less susceptible to unforeseen changes in the completion time of transfers. 

\item Our model takes into account the amount of lateness and uses a nonlinear late delivery penalty.




\item Contrary to commercial freight systems, most of the humanitarian disaster relief organizations work with a certain budget and their goal is to achieve the best possible outcome by spending that budget effectively. Our framework is built based on this practicality.

\item Our framework provides the decision makers with the amount of aid that the system can deliver reliably in certain time frames. This information is useful when coordinating with other organizations involved in disaster relief operations.

\end{enumerate}

\medskip
\singlespacing
\normalsize
\bibliography{references}

\end{document}